%% file: PaperPhDSajo2017.tex
\pdfoutput=1
\documentclass{elsarticle}
\usepackage{amsmath,xcolor,amssymb,amsfonts,multicol,booktabs} %,showkeys
\usepackage[hmargin=2cm,vmargin=2.5cm]{geometry}
\usepackage[english]{babel}
\usepackage[utf8]{inputenc}
\usepackage[plain]{algorithm}
\usepackage{algorithmic}
\newtheorem{theorem}{Theorem}
\newcommand{\set}[1]{\if#1\Omega{\Omega}\else{\mathbb{#1}}\fi}
\let\c\relax
\DeclareMathOperator\c{\mathrm{c}}
\DeclareMathOperator\f{\mathrm{f}}
\DeclareMathOperator\F{\mathrm{F}}
\DeclareMathOperator\h{\mathrm{h}}
\DeclareMathOperator\g{\mathrm{g}}

\DeclareMathOperator\C{\mathrm{C}}
\DeclareMathOperator\V{\mathrm{V}}

\DeclareMathOperator\T{\mathsf{T}}
\DeclareMathOperator\lagrangian{\mathcal{L}}
\DeclareMathOperator\lagrangiana{\mathrm{L}}
\newcommand{\Define}{\mathbin{\raisebox{0.3pt}{:}\!=}}
\DeclareMathOperator{\nnz}{\mathrm{nzz}}
\DeclareMathOperator\diag{{\mathrm{diag}}}

\begin{document}
\title{Preconditioning ideas for the Augmented Lagrangian method}

\author[USB]{A.~M.~Sajo-Castelli}
\ead{asajo@usb.ve}
\address[USB]{Departamento de
C\'omputo Cient\'{\i}fico y Estad\'{\i}stica, Universidad Sim\'on
Bol\'{\i}var, Ap. 89000, Caracas 1080-A,  Venezuela.}

\date{December 16, 2016}

\begin{abstract}
A preconditioning strategy for the Powell-Hestenes-Rockafellar 
Augmented Lagrangian method (ALM) is presented. The scheme exploits the
structure of the Augmented Lagrangian Hessian. It is a modular
preconditioner consisting of two blocks. The first one is associated with the
Lagrangian of the objective while the second administers the Jacobian of the
constraints and possible low-rank corrections to the Hessian.
The proposed updating strategies 
take advantage of
ALM convergence results and avoid frequent refreshing. Constraint 
administration takes into account complementarity over the Lagrange multipliers 
and admits 
relaxation.
The preconditioner is designed for problems where constraint quantity is 
small compared to the search space. A virtue of the scheme is that it is  
agnostic to the preconditioning technique used for the Hessian of the 
Lagrangian function.
The strategy described can be used for linear and nonlinear preconditioning.
Numerical experiments report on spectral properties of  
preconditioned matrices from Matrix Market while some optimization 
problems where taken from the CUTEst collection. Preliminary results indicate 
that the 
proposed scheme could be attractive and further experimentation is encouraged.
\end{abstract}


\maketitle

\section{Introduction}
Augmented Lagrangian methods (ALM) are practical and affordable algorithms 
extensively used in applied fields. They are designed to solve large-scale 
nonlinear optimization problems possibly with nonlinear constraints. 
The constraints 
are classified in two groups: hard and soft. Hard constraints are such that
strict fulfillment is required in order to accept a solution while  
 minor infeasibility for the soft constraints is granted.
Birgin and Martínez in~\cite{libroMLA} present a recent and detailed overview
of practical ALM. These methods are considered a general optimization 
machinery~\cite{Diniz-Ehrhardt2004Augmented,andreani2007augmented,
andreani2008augmented,rockafellar1974augmented,natasa2000,
libroMLA,conn1991globally} in the sense that they successfully cope 
with a great variety of real-life problems. ALM has also been adapted 
to very specific 
applications~\cite{simo1992augmented,fortin2000augmented,afonso2011augmented}.
A great virtue of the method is that it can be accelerated via preconditioning
techniques.

Iterative Krylov-type methods have been used inside ALM and even though they are
very well studied and can handle very large-scale problems, 
it is well known the poor convergence speed. 
Accelerating these methods by preconditioning strategies is common practice.
%
% In its most basic form, preconditioning is considered a subtle mixture of 
% science and art that seeks to transform a complex linear system in such a way 
% to have the same solution but is much more attractive to solve using iterative 
% methods.
%
In the last 50 years a considerable amount of effort has been invested in the
design and construction of effective preconditioners. 
In the context of ALM, there has not been much motivation to study 
acceleration strategies that benefit or exploit convergence results, although
recently some studies on very special problems present non-induced 
preconditioning that benefit from ALM 
convergence~\cite{benzi2011modified,heister2013efficient,benzi2011field}.

State of the art ALM implementations that can be highlighted are 
Algencan~\cite{libroMLA,
andreani2007augmented,andreani2008augmented} and 
Lancelot~B~\cite{gould2003galahad,conn2013lancelot}. 
Both are considered production-grade codes that use a variety of direct and 
iterative methods, such as the Conjugate Gradients method. 
In regards to acceleration, 
Algencan does not offer enough
flexibility while Lancelot~B incorporates a list facilities, it can also 
leave to the user the task of administering the whole preconditioning process.
Surely this last option covers all possible scenarios, but it can also be 
daunting for non-expert users. 
From the user-land perspective, these MLA implementations 
leave a certain void preconditioning-wise. 
The proposed acceleration scheme tries to fill this gap.

In this work we present a modular acceleration scheme that exploits the 
special structure of the Augmented Lagrangian of Powell-Hestenes-Rockafellar.
A key aspect is that the update strategies
take advantage of ALM convergence. 
The scheme has two main ingredients, an \emph{auxiliary} preconditioner 
associated with the Lagrangian function and a machinery that administers the 
constraints.
Separating these components has the great advantages of freedom in choosing the
auxiliary preconditioner and having full control over the Jacobian matrix of 
the constraints and possible low-rank corrections. Such corrections are 
attractive because they promote quality in the approximation to the Hessian of 
the Lagrangian function. The proposed scheme is considered a generalization of
the \textsc{qncgna} preconditioner used in 
Algencan~\cite{birgin2008structured,libroMLA}.

The rest of this document is organized as follows.
We briefly introduce the problem of interest and the Augmented Lagrangian 
Method followed by the presentation of the acceleration scheme. 
The preconditioner is introduced in two parts.
We start by showing how to accelerate ALM for problems with a single constraint,
then a general preconditioner is presented. Update strategies for the various 
components are discussed and some illustrative numerical results are reported.
The work concludes with some final remarks.

\section{The Augmented Lagrangian Method}
Let us consider the following optimization problem
\begin{align}\begin{split}
 \text{minimize}\quad&\f(x)\\
 \text{subject to} \quad&\h(x) = 0,\\
 & \g(x) \le 0,\\
 &x\in\set{\Omega},
 \end{split}\label{eq:problema}
\end{align}
where $\h:\set{R}^{n}\rightarrow\set{R}^{m}$, $\g:\set{R}^{n}\rightarrow\set{R}^{p}$,
$\f:\set{R}^{n}\rightarrow\set{R}$ are continuous functions, in particular
$\f()$ is two times continuously differentiable and
$\set{\Omega}=\{x\in\set{R}^{n} \,|\, \ell\le x\le u\}$.
With the goal to solve~\eqref{eq:problema}, two functions are associated,
the Langrangian
\begin{align}\label{eq:lagrangiano}
 \lagrangian(x, \lambda, \mu) & = \f(x) + \sum_{i=1}^m \h_i(x)\lambda_i  + 
 \sum_{i=1}^p \g_i(x)\mu_i = \f(x) + \h(x)^{\T}\lambda + \g(x)^{\T}\mu,
\end{align}
where $\lambda\in\set{R}^m$ y $\mu\in\set{R}_+^p$ are the 
\emph{Lagrange multipliers} of the problem, and the 
\emph{Augmented Lagrangian} of 
Powell-Hestenes-Rockafellar~\cite{hestenes1969multiplier,rockafellar1973multiplier}
\begin{align}\label{eq:flagaug}
 \lagrangiana_{\rho}(x, \lambda, \mu) &= \f(x) 
 + \frac{\rho}{2} \sum_{i=1}^m \left[\h_i(x) + \frac{\lambda_i}{\rho}\right]^2
 + \frac{\rho}{2} \sum_{i=1}^p \left[\max\bigg(0,\g_i(x) + \frac{\mu_i}{\rho}\bigg)\right]^2\!\!,
\end{align}
with the \emph{external penalty parameter} $\rho>0$, 
$\lambda\in\set{R}^m$ y $\mu\in\set{R}_+^p$. 
Taking the Augmented Lagrangian function as the objective, then 
solving~\eqref{eq:problema} and minimizing~\eqref{eq:flagaug} with respect to
$x\in\set{\Omega}$, are equivalent problems.
The ALM consists of solving a sequence of \emph{sub-problems} and updating the 
the Lagrange multipliers and the external penalty parameter as required.
This naturally divides the iterations in two groups, external or Lagrangian 
iterations and internal iterations.
On the external iterations the multipliers and the penalty parameters are updated while
the internal iterations are dedicated in solving the sub-problem.
Soft constraints are \emph{moved up} to the objective penalizing a shifted
version of infeasibility measure. Hard constraints are enforce inside the 
inner solver which is dedicated to the sub-problem
\begin{align}\begin{split}\label{eq:minlagaug}
 \text{minimize}\quad & \lagrangiana_{\rho_k}(x, \lambda_k, \mu_k)\\
 \text{subject to} \quad & x\in\set{\Omega},
 \end{split}
\end{align}
where $\rho_k$, $\lambda_k$ and $\mu_k$ are fixed.

The ALM is conceptually presented in Algorithm~\ref{alg:MLA}.
The algorithm is schematic and leaves open the choice on 
how to solve the sub-problem~\eqref{eq:minlagaug}, it only requires
that $x_k$ be an \emph{approximate} solution.
On each external iteration (Step~4) it is considered if $x_k$ 
has done enough progress in regards to feasibility and complementarity.
In cases where progress is good, the external penalty parameter 
does not need to change, on the contrary the parameter must be incremented.
The idea behind the shifts ${\lambda}/\rho$ and ${\mu}/\rho$ is
that even when the external penalty parameter has a \emph{moderate}
penalty, there can exist \emph{adequate} values for the shifts on the 
multipliers where it is possible to find an acceptable solution 
to~\eqref{eq:minlagaug}.
Given the interest of the algorithm to produce a numerically 
attractive result, on Step~5 the multipliers are safeguarded.

\input{alg-MLA}

We briefly present convergence results for Algorithm~\ref{alg:MLA}, details and 
theorem proofs are given in~\cite[Chapter~5]{libroMLA}.
For the following results it is assumed that on Step~1 of 
Algorithm~\ref{alg:MLA}, $x_k$ is a global minimizer 
of~\eqref{eq:minlagaug}:
\[ \lagrangiana_{\rho^k}(x^k, \bar{\lambda}^k, \bar{\mu}^k) \leq 
\lagrangiana_{\rho^k}(x, \bar{\lambda}^k, \bar{\mu}^k) + \varepsilon^k, 
\quad \forall\, x\in\Omega,\]
where $\{ \varepsilon^k\}\subseteq\set{R}_+$ is bounded and 
$\varepsilon^k$ need not be \emph{small}.

\begin{theorem}[Feasibility]
Let $\{x^k\}$ be a sequence generated by Algorithm~\ref{alg:MLA} under
the previous assumption and let $x^*$ be a limit point of this sequence.
Then,
\[\|h(x^*)\|^2_2 + \|g(x^*)_+\|^2_2 \leq \|h(x)\|^2_2 + \|g(x)_+\|^2_2, \quad \forall\, x\in\Omega.\]
\end{theorem}
This result guarantees that if the problem is feasible, then every 
limit point of the sequence generated by ALM is also feasible.

\begin{theorem}[Optimality]
Let $\{x^k\}$ be a sequence generated by Algorithm~\ref{alg:MLA} and 
$x^*$ be a limit point of $\{x^k\}$.
Suppose that $\varepsilon^k\xrightarrow{k\rightarrow\infty}0$ and that
the problem~\eqref{eq:minlagaug} is feasible. On Step~4, when possible, 
$\rho_{k+1}$ is not incremented. Then, $x^*$ is a global minimizer.
\end{theorem}

For practical purposes it is not necessary to differentiate between the equality
and inequality constraints, they can be both considered under the 
\emph{umbrella} constraint function $\c(x)$ and the associated multiplier 
vector $\lambda$. Optimality conditions and constraint handling requires only
minor attention. The Augmented Lagrangian becomes
\begin{align}
\begin{split}
 \F(x)=\lagrangiana_{\rho}(x, \lambda) &= \f(x) + \frac{\rho}{2} \sum_{i\in\set{E}} \left[c_i(x) + \frac{\lambda_{i}}{\rho}\right]^2
 + \frac{\rho}{2} \sum_{i\in\set{I}} \left[\max\bigg(0,c_i(x) + \frac{\lambda_{i}}{\rho}\bigg)\right]^2
\end{split}\label{eq:lagaug}
\end{align}
where the sets $\set{E}$ and $\set{I}$ contain the indexes 
of equality and inequality constraints respectively.

Supposing that problem~\eqref{eq:problema} has solution, then
the pair $(x_k, \lambda_k)$ is minimizer of~\eqref{eq:lagaug} 
---and in consequence is also of~\eqref{eq:flagaug} and solution
to~\eqref{eq:problema}--- when it satisfies
\begin{align}
 \left\| 
 P_{\set{\Omega}}\left(x_k-\left[
\nabla\lagrangiana_{\rho_k}(x_k, \lambda_k) + \sum_{j=1}^m \lambda_k\,\nabla\c_j(x_k)
\right]\right) -x_k \right\| & \le\varepsilon_{\text{opt}},\label{eq:KKT1}\\
\max\bigg(\max_{j\in\set{E}}\Big(\big|c_j(x_k)\big|, \max_{j\in\set{I}} 
\Big(\big|\min(-c_j(x_k), \lambda_k)\big|\Big)\bigg) &
\le\varepsilon_{\text{fact}},\label{eq:KKT2}\\
\max\bigg(\max_{j\in\set{E}}\Big(\big|c_j(x_k)\big|, \max_{j\in\set{I}} 
\Big(\c_j(x_k)_+\Big)\bigg) &\le\varepsilon_{\text{fact}}\label{eq:KKT3},
\end{align}
where
\begin{align*}
\lambda_k=\begin{cases}
           \lambda_{j,k-1}+\rho_k\,\c_j(x_k) & \quad\text{if } j\in\set{E} 
           \text{ or } \lambda_{j,k-1}+\rho_k\,\c_j(x_k)>0,\\
           0 & \quad\text{otherwise.}
          \end{cases}
\end{align*}
The operator $P_{\set{\Omega}}$ is the Euclidean projection on the convex set
$\set{\Omega}$. Conditions~\eqref{eq:KKT1} and~\eqref{eq:KKT2}
insure that $(x_k, \lambda_k)$ is a stationary point of the problem while
condition~\eqref{eq:KKT3} guarantees that the point satisfies the 
required feasibility tolerance.

It is advantageous to differentiate between Lagrangian and internal 
iterations.
External iterates and counter are identified by $x_k$ and $k$ respectively while
internal ones use $z_{\ell}$ and $\ell$ respectively.
If the set $\set{\Omega}=\set{R}^n$ then the sub-problem is
said to be \emph{unconstrained}, on the contrary it will be assumed to be 
\emph{convex constrained}.

Practical Newton-type methods are very popular for solving unconstrained 
problems and require a descent direction of first order 
given by solving the quadratic model
\begin{align*}
 H_{\ell}\, d_{\ell} &= -\nabla\F(z_{\ell}),
\end{align*}
where $H_{\ell}$ is the Hessian $\nabla^2\F(z_{\ell})$ or 
an approximation to it.
For the Truncated-Newton method,  $d_{\ell}$ is estimated using the 
Conjugate Gradients method in the SPD case or Minimal Residual method when
$H_{\ell}$ is symmetric but undefined. These methods are quite attractive if
accelerated with high quality preconditioners.

Practical solvers for convex constrained problems are the 
Spectral Projected Gradient method 
(SPG)~\cite{birgin2000nonmonotone,birgin2001algorithm,
birgin2014spectral} and its preconditioned variant (PSPG)~\cite{bello2005precond}.
PSPG has received some attention~\cite{birgin2013PSPG} and 
can be viewed as a nonlinear preconditioned variant that combines the Preconditioned 
Spectral Gradient method and SPG.
For these type of problems we assume that $P_{\Omega}$ is the Euclidean projection operator 
over the convex set $\set{\Omega}$, exists and is of acceptable cost. We also require
that the first order derivatives of $\f()$ and $\c()$ exist wherever required.

\section{Preconditioning ideas}
Under our context, accelerating ALM may refer to accelerate the resolution of the
sub-problem (Newton-type directions) as well as the acceleration for the 
estimation of the descent direction (Cauchy-type machinery). Given the nature of the
method to solve at each external iteration an optimization problem, and that the
sequence of these problems tend to be similar, preconditioning schemes must 
\emph{recycle} between iterations in order to be attractive.

In order to propose preconditioning schemes for \eqref{eq:minlagaug}, it is 
necessary to study the explicit form of the Hessian 
of the objective (Augmented Lagrangian function).
For historical reasons, most applications use and implement the Lagrangian 
function and not its augmented counterpart. Although, as we shall see both are
closely related.
The Lagrangian function, gradient and Hessian associated to~\eqref{eq:lagaug}, are
\begin{align*}
 \lagrangian(z, \lambda) & = \f(z) + \sum_{i=1}^m \c_i(z)\lambda_i  = \f(z) + \c(z)^{\T}\lambda,\\ % \quad \text {and}\\
 \nabla\lagrangian(z, \lambda) & 
= \nabla\f(z) + \sum_{i=1}^m \lambda_i\nabla\c_i(z)=
\nabla\f(z) + \nabla\c(z)\lambda, \quad \text{with} & 
 \nabla \c(z) &= 
\begin{bmatrix}
| & | & & |\\
\nabla \c_1(z) & \nabla \c_2(z) & \cdots & \nabla \c_m(z)\\
| & | & & |
\end{bmatrix},\\
\nabla^2\lagrangian(z, \lambda) & = \nabla^2\f(z) + \sum_{i=1}^m \lambda_i\nabla^2\c_i(z).
\end{align*}
On the other hand, the Augmented Lagrangian counterpart is
\begin{align*}
 \lagrangiana_{\rho}(z, \lambda) &= \f(z) + \frac{\rho}{2} \sum_{i=1}^m \left[c_i(z) + \frac{\lambda_i}{\rho}\right]^2
 =\f(z) + \c(z)^{\T}\lambda + \frac{\rho}{2} \c(z)^{\T}\c(z) + \frac{1}{2\rho}\lambda^{\T}\lambda,\\
 \nabla \lagrangiana_{\rho}(z, \lambda) % &= \nabla\f(z) + \sum_{i=1}^m \lambda_i\nabla\c_i(z)  + \rho \sum_{i=1}^m \c_i(z) \nabla\c_i(z)\\
 &=\nabla\f(z) + \sum_{i=1}^m \big[\lambda_i  + \rho \c_i(z) \big] \nabla\c_i(z)
=\nabla\f(z) + \nabla\c(z)\,\big[\lambda  + \rho \c(z) \big],\\
 \nabla^2 \lagrangiana_{\rho}(z, \lambda) &= \nabla^2\f(z) + \sum_{i=1}^m \lambda_i\nabla^2\c_i(z) + 
 \rho \sum_{i=1}^m \c_i(z) \nabla^2\c_i(z) + \rho \sum_{i=1}^m \nabla \c_i(z)\nabla\c_i(z)^{\T},
 \intertext{grouping $\nabla^2\c_i(z)$ %}
 and doing the change: $\hat{\lambda}=\lambda+\rho\c(z)$, we have}
 \nabla \lagrangiana_{\rho}(z, \lambda) &= \nabla^2\f(z) + \sum_{i=1}^m \hat{\lambda}_i \nabla^2\c_i(z) + \rho \nabla\c(z) \nabla\c(z)^{\T}\!.
\end{align*}
Noting that $\nabla^2\lagrangian(z, \hat{\lambda}) = \nabla^2\f(z) + \sum_{i=1}^m \hat{\lambda}_i \nabla^2\c_i(z)$, we obtain the key identity
 \begin{align*}
\nabla^2 \lagrangiana_{\rho}(z, \lambda) &= \nabla^2\lagrangian(z, \hat{\lambda}) +  \rho \nabla\c(z) \nabla\c(z)^{\T}, \quad \text{or equivalently,}\quad \nabla^2\F(z_{\ell}) = H_{\ell} = M_{\ell} + \rho\,V_{\ell}{V_{\ell}}^{\T}.
\end{align*}
This matrix sum has many relevant characteristics. For instance, since $M$ is the
sum of Hessian matrices it is symmetric and close to a solution is 
definite~\cite{nocedal2006numerical}. Under certain choices of $M$, 
 it can be verified that it is always definite~\cite{natasa2000}. 
The Gauss-Newton matrix $VV^{\T}$  is always 
symmetric and can be regarded 
semi-definite. Moreover, if $m<n$ then the rank of $VV^{\T}$ is at best $m$,
the number of soft constraints. In practice, methods exploit the complementary
condition in such a way that the rank of $VV^{\T}$ is at most the number of
\emph{active non-relaxed} constraint count at current iterate. 

When solving problem~\eqref{eq:minlagaug} using Newton-type directions, it is
necessary to solve the quadratic model
\begin{align}
\nabla^2\F(z)\,d = H d = \big[M + \rho VV^{\T}\big]\,d = -\nabla\F(z).\label{eq:LSD}
\end{align}
For Cauchy-type directions, if $P$ is a preconditioner for $H$, then we 
interpret the preconditioner $P$ as an 
approximation to $H^{-1}$ and the \emph{enriched} descent direction is 
\begin{align*}
 d &= -P\,\nabla\F(z).
\end{align*}
Contrasting the previous two expressions,
applying nonlinear preconditioning for Gradient-type methods
fortunately is analogous as accelerating the Newton-type machinery.
Our interest lies in solving the linear system~\eqref{eq:LSD} using Krylov-type iterative methods,
PCG when $H$ is PD or MinRes for the undefined case.
Preconditioning schemes for ALM must be able to at least exploit the following
two desirable key features:
\begin{enumerate}
 \item \emph{Preconditioner recycling.} The idea is to take advantage of 
 convergence for the sub-problems. Supposing that \eqref{eq:problema} has 
 solution, then we expect
 $H_k\xrightarrow{k\rightarrow \infty}{H^*}$,
 meaning that starting from a certain iterate
 $k$ (or $\ell$) it is possible to bound the difference
 $H_{k+1} - H_{k}$ and successfully use the same preconditioner onwards.
 \item \emph{Preconditioner update.} Assembly of $P_{k+1}$ should partially 
 reuse work invested in assembling $P_{k}$. Generally speaking this is not 
 straightforward. As an illustrative example on the involved difficulties,
 low-rank 
 updates~\cite{bellavia15,birgin2008structured,Mar95,Mar93,conn1991convergence}
 are highlighted
 with special focus on BFGS-type 
 corrections~\cite{nocedal1980updating,andrei2007scaled}.
\end{enumerate}

In what follows we propose a new preconditioning scheme that takes advantage of
these features.

\subsection{Preconditioning singly constrained sub-problems}\label{sec:precond1}

Many applications can be modeled using singly constrained problems: support
vector machine formulations for classification or pattern 
recognition~\cite{cores2009on,dai2006new} are but two mainstream examples.
This subsection introduces a new inverse preconditioner scheme for solving 
the following singly constrained problem
\begin{align}
\begin{split}
 \text{minimize}\quad&
  \F(z)=\lagrangiana_{\rho}(z, \lambda_1) = \f(z) + \frac{\rho}{2} 
  \left[c_1(z) + \frac{\lambda_{1}}{\rho}\right]^2\\
 \text{subject to}\quad&z\in\set{\Omega}.
 \end{split}\label{eq:minlagaugc1}
\end{align}
We simplify notation by relaxing
the index from the constraint and its associated Lagrange multiplier.
With the idea to introduce the preconditioner, suppose a Newton-type
machinery is used to solve~\eqref{eq:minlagaugc1}, then 
at each internal iteration $\ell$ it is required to find $d$ from the linear system
\[
H_{\ell}\;d = -\nabla\F(z_{\ell}),
\]with\[
H_{\ell}=\nabla^2\lagrangiana_{\rho_k}(z_{\ell}, \lambda_k) = \nabla^2\lagrangian\big(z_{\ell}, \lambda_k+\rho_k\c(z_\ell)\big) + \rho_k\; \nabla\c(z_{\ell})\,\nabla\c(z_{\ell})^{\T} 
= M_{\ell}+\rho_k\, v_{\ell}{v_{\ell}}^{\T}.
\]
Relaxing the iteration indexes, we have
\begin{align}
 \big[M+\rho\,vv^{\T}\big]\; d = -\nabla\F(z),\label{eq:LSD1}
\end{align}
where $v=\nabla\kern-1.25pt\c(z)$ is considered a non-null vector, $vv^{\T}$ is
a rank-1 matrix and we suppose $M$ is full rank.
Given our interest in solving~\eqref{eq:LSD1} using PCG or MinRes, applying the
preconditioner~$P^{-1}$~\cite[Algorithm~5.3: line~3 and equation~(5.38d)]{nocedal2006numerical} 
reduces to efficiently compute the product
\[ h_j = P^{-1} r_j, \quad \text{with} \quad r_0=H\,z_0+\nabla\F(z_0).\] 
The matrix in~\eqref{eq:LSD1} is special in the sense that it is the sum of an
invertible matrix plus a rank-1 matrix.
Suppose that $P^{-1}$ is precisely 
$H^{-1} = \Big[M+\rho\,vv^{\T}\Big]^{-1}\!\!\!\!\!,\;\;$ 
then the product $h_j = P^{-1} r_j$ can be computed using the
Sherman-Morrison identity~\cite{sherman1950}, 
\begin{align*}
 h_j &= \Big[M+\rho\,vv^{\T}\Big]^{-1} r_j,\\
     &= \bigg[M^{-1} - \frac{1}{1+\rho\,v^{\T}M^{-1}v} \;\rho\,M^{-1}vv^{\T}M^{-1} \bigg] r_j = M^{-1}r_j - \frac{1}{1+\rho\,v^{\T}M^{-1}v} \;\rho\,M^{-1}vv^{\T}M^{-1}  r_j,
\end{align*}
solving for $a$, {$Ma=r_j$} or {$a=M^{-1}r_j$}, we obtain
\[ h_j = a -\frac{\rho\,}{1+\rho\,v^{\T}M^{-1}v} \;M^{-1}vv^{\T}a,\]
analogously we find  $b$, {$Mb=v$} or {$b=M^{-1}v$} 
and observing that $v^{\T}a$ is a scalar, we finally have
\[h_j = a -\frac{\rho\,v^{\T}a}{1+\rho\,v^{\T}b} \;b.\]
If $M$ is not trivially invertible, finding $a$ and $b$ requires an 
\emph{auxiliary} preconditioner $P_M$. Each form of obtaining $a$ and 
consequently $b$, give rise the different preconditioning strategies.
It is noteworthy to remark that the preconditioner $P$ is never assembled and is
considered an \emph{abstract} preconditioner which relies upon the auxiliary
preconditioner $P_M$. This highlights the \emph{agnostic} nature of $P$,
for it does not enforce any specific choice over $P_M$.
Under this scheme, the spectrum of the preconditioned 
matrix $P^{-1}\big[M+\rho\,vv^{\T}\big]$ is described by the following result.
\begin{theorem}
Suppose $M$ is the Hessian matrix of the Lagrangian function associated 
to~\eqref{eq:minlagaugc1}, let the preconditioner for $H$
be $P=\big[M+\rho\,vv^{\T}\big]$, furthermore let $P_M^{-1}>0$ be a preconditioner for
$M$, $\rho\ge1$ and $\mathcal{E}_M \Define P_M^{-1}M - I$, then the spectrum of $P^{-1}H$ is
$$\Lambda\big(P^{-1}H\big) = \Lambda\Big(I + \big(1 -  \upsilon \big) \;P_M^{-1} vv^{\T}\mathcal{E}_M\Big),\quad \text{ where }\quad
\upsilon = \frac{\rho }{1 + \rho v^{\T}P_M^{-1}v}.$$
\end{theorem}
The proof of this theorem can be found in~\cite{thesis-phd-usb-sajo-2016}.
This result has a few important implications. If 
$\rho \rightarrow \infty$, then $P^{-1}H \rightarrow I.$ This indicates
that for very large values of $\rho$ it is not attractive or even necessary to
precondition. It also shows that the condition of $P^{-1}H$ can be described
in terms of $P_M^{-1}M$, this is to say that the quality of $P$ is given
in direct relation to the quality of $P_M$.

\subsection{General Preconditioner}\label{sec:precondm}\noindent
In this section we work with the problem
\begin{align}
\begin{split}
 \text{minimize}\quad&
  \F(z)=\lagrangiana_{\rho}(z, \lambda) = \f(z) + \frac{\rho}{2} 
  \sum_{i=1}^m \left[c_i(z) + \frac{\lambda_{i}}{\rho}\right]^2\\
 \text{subject to}\quad&z\in\set{\Omega},
 \end{split}\label{eq:minlagaugcm}
\end{align}
where the quantity of constraints is less than the problem dimension.
Analogous to~\eqref{eq:LSD1}, the linear system to solve is
\[H d = \big[M_{}+\rho\,VV^{\T}\big]\, d = -\nabla\F(z),\]
where $V=\nabla\c(z)$ is the constraint Jacobian, 
 $VV^{\T}$ is the Gauss-Newton matrix of rank at most $m$ and we suppose 
 that $M$ is rank complete.
Initially, we set the preconditioner to $P^{-1}=\Big[M+\rho VV^{\T}\Big]^{-1}$.
The use of the 
Sherman-Morrison-Woodbury~\cite{woodbury1950inverting,hager1989updating} 
inverse closed formula rises many difficulties, specially
the stringent requirement of a composed matrix inverse which for our case is 
not practical.
Fortunately, $H$ can be rewritten as a recursion over the columns of $V$, suppose~\eqref{eq:minlagaugcm} has three constraints ($m=3$), then
 \[ H = \Big(\big((M+\rho v_1v_1^{\T})\;\;+\;\;\rho v_2v_2^{\T}\big)\;\;+\;\;\rho v_3v_3^{\T}\Big).\]
This formulation shows that $H$ can always be written as an invertible matrix
plus a rank-1 matrix. This is also true for the preconditioner $P$,
let $P_0\Define M$, then
\begin{align*}
 P_1 &= P_0 + \rho\,v_1v_1^{\T},&
 P_2 &= P_1 + \rho\,v_2v_2^{\T},&
 \text{and finally,}&&
 P = P_3 &= P_2 + \rho\,v_3v_3^{\T},
\end{align*}
which leads to the following recursion,
$$P = P_m=P_{m-1}+\rho\,v_mv_m^{\T}, \text{ with } P_0:=M,\;\text{and } m\ge1.$$ 
Using Miller's inverse formula~\cite{miller1981inverse}, the product 
$h_j=P^{-1}r_j$ is reduced to the recursion
\begin{align}
 h_j=P^{-1}r_j &= P_3^{-1}r_j = \big[P_2 + v_3{v_3}^{\T}\big]^{-1}r_j,\label{eq:hjrec}
 \end{align}
but this was previously shown how to be solved, and is 
 \begin{align*}
        h_j=a_4&= a_3 -\frac{\rho\,v_3^{\T}a_3}{1+\rho\,v_3^{\T}b_3} \;b_3, \text{ with } a_3=P_2^{-1}r_j \text{ and } b_3 = P_2^{-1}v_3.
\intertext{Now, in a similar fashion we compute $a_i,\; i=3:1$:}
 a_3&=P_2^{-1}r_j = a_2 -\frac{\rho\,v_2^{\T}a_2}{1+\rho\,v_2^{\T}b_2} \;b_2, \text{ with } a_2=P_1^{-1}r_j \text{ and } b_2 = P_1^{-1}v_2,\\
 a_2&=P_1^{-1}r_j = a_1 -\frac{\rho\,v_1^{\T}a_1}{1+\rho\,v_1^{\T}b_1} \;b_1, \text{ with } a_1=P_0^{-1}r_j \text{ and }  b_1 = P_0^{-1}v_1,\\
 a_1&=M^{-1}r_j.
 \intertext{Unfortunately, the estimation of the $b_i, \;i=3:1$ is a bit more elaborate: in order to compute $b_3$ it is required to have  $P_2^{-1}v_3$, which in turn
 requires  $P_1^{-1}v_3$ and $P_1^{-1}v_2$,}
b_3 &= P_2^{-1}v_3 = a'_2 -\frac{\rho\,v_2^{\T}a'_2}{1+\rho\,v_2^{\T}b_2} \;b_2, \text{ with } a'_2=P_1^{-1}v_3 \text{ and } b_2 = P_1^{-1}v_2,\\
b_2 &= P_1^{-1}v_2 = a'_1 -\frac{\rho\,v_1^{\T}a'_1}{1+\rho\,v_1^{\T}b_1} \;b_1, \text{ with } a'_1=P_0^{-1}v_2 \text{ and } b_1 = P_0^{-1}v_1,\\
b_1 &= M^{-1}v_1.
\end{align*}
Implementation-wise, the $a'_i$ are computed at the same time as the $b_i,\;i=2:1$ 
which leads to a secondary recursion.

Although the previous formulation requires double recursion over the columns of 
$V$, it will be shown to have many attractive features.
The elements $a_i, b_i$, $a'_i$ and finally the acceleration product $h_j$, 
are computed by applying multiple times the Sherman-Morrison identity. This is
considered a very attractive aspect because it only requires Matrix-Vector and
internal products.

In the same spirit as the single constrained case, each form of obtaining
$a_i$ and $b_i$ give rise to different preconditioners and this variant is also
considered to be an agnostic acceleration scheme which leaves open the 
strategy on how to estimate $M^{-1}{a_i}$ and $M^{-1}{b_i}$.
This has numerous advantages, we highlight the fact that handling the Lagrangian 
independently from the constraint Jacobian permits to exploit the sparse 
structure of the problem.

In the next section we generalize the previous example of three constraints.

\subsubsection{The $B$ Matrix}\label{sec:B}\label{sec:evalsmm}
Observing the recursion that computes the product~\eqref{eq:hjrec}, it can be 
seen that elements $a'_i$ and $b_i$ \emph{do not} depend on 
$r_j$ but on the columns of $V$. This in a  natural way induces to 
\emph{pre-compute} $a'_i$ and $b_i$, and save them efficiently in a $B$
\emph{storage matrix}. In what follows we show how to efficiently compute
the product $h=P^{-1}r$.

Let $h_0 \Define M^{-1}r$, $B_{0;*}=M^{-1}V$ where $B_{i;j}$ indicates the $j$-th column of matrix $B_i$.

\noindent 
On a first pass, storage matrix $B$ is assembled,
\begin{align*}
B_{i;j} &= B_{i-1;j} - \frac{\rho v_i^{\T} B_{i-1;j}}{1 + \rho\,v_i^{\T}B_{i-1;i}}B_{i-1;i}, \quad \text{ for } j=i:m,\; 
\text{ with } i = 1:m.
\intertext{On a second pass, $h$ is computed,}
 h_i &= h_{i-1} - \frac{\rho\,v_i^{\T} h_{i-1} }{ 1 + \rho\,v_i^{\T} B_{i-1;i}} B_{i-1;i}, \quad \text{ for } i=1:m.
\end{align*}
Finally $h= h_m$.
It is important to note that from the implementation point of view, the elements
of $B_{i+1;*}$ overwrite those of $B_{i;*}$ in such a way as not to waist 
 storage. The recursion to compute $h$ distills to
\[h_i = h_{i-1} - \frac{\rho\,v_i^{\T} h_{i-1} }{ 1 + \rho\,v_i^{\T} B_{i}} B_{i},\]
where $B_i$ regains its classical meaning indicating the  $i$-th column of $B$.
Storage matrix $B$ unites in a single matrix all the ingredients to 
estimate $h$ and as such it is intimately related to $M$ and $V$.
The previous formulation efficiently computes the product
$P^{-1}w$ or $[M+\rho\,VV^{\T}]^{-1}w$ for any $w$ as long as $M$ and $V$
stay relatively the same.
Significant changes in $V$ and/or $M$ force re-assembly of $B$.
This fact outlines certain updating aspects that an acceleration scheme must 
be aware of.
It is evident that this strategy to compute $h$ is only attractive when
$m$ is far from $n$. 
The proposed preconditioner $P$ can be shown trivially that is adequate for 
 use within PCG and MinRes.

\subsubsection{Secant type directions}
In large-scale applications computing the Augmented Lagrangian Hessian is a 
luxury seldom available, in general terms a reasonable approximation is used.
Within our context, the Augmented Lagrangian Hessian has clear differentiation
between its components
\begin{align*}
 \nabla^2 \lagrangiana_{\rho}(z, \lambda) &= \nabla^2\f(z) + \sum_{i=1}^m [\lambda_i + \rho\,\c_i(z)]\,\nabla^2\c_i(z)
 + \rho \sum_{i=1}^m \nabla \c_i(z)\nabla\c_i(z)^{\T}
 =\nabla^2\f(z)+\C(z)+\V(z).
\end{align*}
Kreji\'{c} et al.%
% and colleagues
~\cite{natasa2000} propose using $\nabla^2\f(z)+\V(z)$ and
Birgin and Martínez suggest in~\cite{libroMLA} to use only $\V(z)$ along with 
two \emph{corrections}.
These involve a spectral correction~\cite{raydan1993barzilai,
Raydan1997SC2,birgin2000nonmonotone,
fletcher1990low} using the associated Rayleigh quotient~\cite{horn1985matrix,trefethen1997numerical} and 
a second correction~\cite{birgin2008structured,libroMLA} in the spirit of BFGS 
that forces to satisfy the Secant equation~\cite{Mar93,Mar95}
\begin{gather*}H_{\ell}s_{\ell}=y_{\ell},\\
 \hat{H}_{\ell}=\nabla^2\f(z_{\ell})+\V(z_{\ell}), \qquad s=z_{\ell}-z_{\ell-1}, \qquad y_{\ell}=\nabla\F(z_{\ell})-\nabla\F(_{\ell-1}).
\end{gather*}
These corrections are low-rank~\cite{conn1991convergence,
babaie2011modified,leong2010monotone,andrei2007scaled,Martinez1990family,Li2001a} 
and are crafted with the main purpose of having a closed inverse form. 
For illustrative purposes, $\hat{H}$ is corrected spectrally and with the famous
BFGS formula,
\begin{align*}
 \hat{H}_+ &\Define \hat{H} + \sigma\,I, \quad \text{ with } \sigma = \frac{(y-\hat{H}(z_{\ell})\,s)^{\T}s}{s^{\T}s},\\
 H &= \hat{H}_+ + \frac{yy^{\T}}{s^{\T}y} - \frac{\hat{H}_+ss^{\T}\hat{H}_+}{s^{\T}\hat{H}_+s}, \quad \text{if $s^{\T}y\ne0$.}
\end{align*}
Now, if $P^{-1} \approx H^{-1}$, the assembly and use of the preconditioner is 
analogous as shown at the start of \S\ref{sec:precondm} since $H$ can be rewritten 
as an invertible matrix plus matrices of rank-1. Let us see,
\begin{align*}
  H &= \hat{H}_+ + \frac{yy^{\T}}{s^{\T}y} - \frac{\hat{H}_+ss^{\T}\hat{H}_+}{s^{\T}\hat{H}_+s}, \quad
  \text{noting that $\hat{H}_+$ is symmetric, we have}\\
   &= \nabla^2\f(z)+\V(z_{\ell}) +\sigma\,I + \nu\,yy^{\T} - \psi\,ww^{\T}, \qquad
   \nu\Define\frac{1}{s^{\T}y},\quad w\Define \hat{H}_+s,\quad \psi\Define\frac{1}{s^{\T}w}
   \\
   &=\nabla^2\f(z)+\sigma\,I + VV^{\T} + \nu\,yy^{\T} - \psi\,ww^{\T}, 
\intertext{letting $M=\nabla^2\f(z)+\sigma\,I$, we have}
H & = M + VV^{\T} + \nu\,yy^{\T} - \psi\,ww^{\T}.
\end{align*}
That is, $H$ has the form required to assemble $P^{-1}$. 
The matrix $M$ is guaranteed to be rank complete by means of the first
spectral correction. With some abuse of notation, the constraint Jacobian
is \emph{augmented} to accommodate the elements associated with the 
 BFGS correction
$$V \leftarrow \left[\begin{array}{ccc}
                       | & | & | \\
                       \sqrt{\rho}\;V & \sqrt{\nu}\;y & \sqrt{\psi}\;w\\
                       | & | & | 
                      \end{array}\right], \quad \text{signs}:= [1, 1, \dots, 1, 1, -1]^{\T}.$$ 
The auxiliary vector ``signs'' is used inside the recursion associated with $B$.

\subsubsection{Special Exact Case}\label{sec:exact}
When using Quasi-Newton directions, the choice of the approximation to the
Lagrangian Hessian plays a key role in determining the quality
of the auxiliary preconditioner associated to $M$.
If the election of $M$ that approximates $\nabla^2\!\lagrangian(z)$ has explicit
inverse, then $P^{-1}$ is the exact (theoretical) inverse of $H$. A trivial 
case is to choose $M$ diagonal.
As a side-effect under this context, the proposed preconditioner can be seen
as a generalization of the \textsc{qncgna} preconditioner of Birgin and 
Martínez~\cite{birgin2008structured,libroMLA}.
The preconditioner is designed to work on the linear system 
\begin{gather*}
V(z_{\ell}) \; d = -g_{\ell}. %\\
% \nabla^2\lagrangiana_{\rho_k}(z, \lambda_k)\approx V(z)=\rho \sum_{i=1}^m \nabla \c_i(z)\nabla\c_i(z)^{\T}.
\end{gather*}
The matrix $V(z)$ is symmetric and SPD, it is corrected spectrally and, if possible, 
a second BFGS-style correction is applied.
The corrected approximation to the Augmented Lagrangian Hessian is
\begin{align}
 H = \left \{\begin{array}{lcl}
      V_+:=V(z) + \sigma\,I && \text{if } y^{T}s < 10^{-8}\,\|y\|\,\|s\|,\\
      V_+ + \frac{yy^{\T}}{s^{\T}y} - \frac{V_+ss^{\T}V_+}{s^{\T}V_+s} && \text{on the contrary,}
\end{array}\right.\label{eq:Hcase}
\end{align}
and the \textsc{qncgna} preconditioner
\begin{align*}
 P_{\text{\textsc{qncgna}}} = \left \{\begin{array}{lcl}
      D_+:=\diag\big(V(z)\big) + \sigma_D\,I && \text{if } y^{T}s < 10^{-8}\,\|y\|\,\|s\|,\\
      D_+ + \frac{yy^{\T}}{s^{\T}y} - \frac{D_+ss^{\T}D_+}{s^{\T}D_+s} && \text{on the contrary.}
\end{array}\right.
\end{align*}
Witch has explicit closed inverse
\[ P_{\text{\textsc{qncgna}}}^{-1} = D^{-1}_+ + \frac{(s- D^{-1}_+y)s^{\T}+s(s- D^{-1}_+y)^{\T}}{s^{\T}y} - \frac{(s- D^{-1}_+y)^{\T}yss^{\T}}{(s^{\T}y)^2}.\]
Note that the search direction is done over $V\!(z)$  and the preconditioner is 
over $\diag\!\big(V\!(z)\big)$.

Now taking up~\eqref{eq:Hcase}, the matrix $H$ can be rewritten
as a rank complete matrix plus the sum of rank-1 matrices. Suppose the BFGS
correction is possible, then
\begin{align*}
H &= V(z) + \sigma\, I  + \frac{yy^{\T}}{s^{\T}y} - \frac{V_+ss^{\T}V_+}{s^{\T}V_+s}, \quad V_+=V(z) + \sigma\,I.
\intertext{Noting that $V_+$ is symmetric, we have}
&= V(z) + \sigma\, I  + \frac{yy^{\T}}{s^{\T}y} - \frac{V_+s(V_+s)^{\T}}{s^{\T}V_+s},\\
  &= \sigma\, I + \rho\,VV^{\T} + \nu yy^{\T} - \psi ww^{\T}, \quad \nu= \frac{1}{s^{\T}y},\quad w=V_+s,\quad \psi=\frac{1}{{s^{\T}w}}.
\end{align*}
Let $M=\sigma\, I$ and $V$ ---again, with some notation abuse--- be the matrix 
of size $n\times(m+q)$ that gathers the constraint Jacobian and the $q$ vectors 
associated with the second BFGS correction, then
defining trivially $P_M = \sigma^{-1}I$ we have that 
$P^{-1}$ is the explicit closed inverse of $H$.

\subsection{Update Strategies}\label{sec:upd}
There exist a variety of updating choices for the 
preconditioner within ALM. Given the available granularity inside the recursive nature of
applying the preconditioner ($P^{-1}r_k$), update 
strategies enjoy a fine-grain control over each component. Practical 
update strategies monitor changes on 
$\big\|M_{\ell}-M_{\ell-1}\big\|$ and $\big\|V_{\ell}-V_{\ell-1}\big\|$ 
independently and take the following actions:
\begin{enumerate}
 \item Update $P_M$ whenever  
 $\big\|M_{\ell}-M_{\ell-1}\big\|_1>\delta_M$.
 \item Update $B$ each time  
 $\big\|V_{\ell}-V_{\ell-1}\big\|_1>\delta_V$ or
 $\big\|M_{\ell}-M_{\ell-1}\big\|_1>\delta_M$.
 \item Apply \emph{relaxation} over the columns of $V$ when  
 assembling 
 $B$ or at the moment of applying the preconditioner ($P^{-1}r_k$). This can be 
 done by 
  using only those columns of $V$ that in norm are greater that a 
  threshold or that have huge infeasibility measure,\label{itm:a}
  $V=(v_i),\; i=1:m\;\;\big|\;\; \|v_i\| > \varepsilon_v>0 \;\vee\;  i\in\set{E} \;\;\big|\;\;|\c_i(x)|>\varepsilon_c>0 \;\vee$ %\\
  $ i\in\set{I} \;\;\big|\;\;\big(\c_i(x)\big)_+>\varepsilon_c$. 
\end{enumerate}
The idea behind the first item is that 
if $M_{\ell}$ and $M_{\ell+1}$ are similar, then possibly 
$P_{M_{\ell}}$ will also be a good preconditioner for $M_{\ell+1}$.
Item two establishes that $B$ should only be updated if 
$V_{\ell+1}$ and $V_{\ell}$ greatly differ. It also forces an update if
$P_{M_{\ell}}$ was updated. Up until not finding the final search space,
$V$ will be changing drastically and update schemes must be aware of this.
The last item is based on the idea that small elements should have small
contributions and can be safely discarded.

\subsubsection{Strategies for $P_M$}\label{sec:updM}
The proposed scheme leaves open the choice on how to precondition $M$.
Clearly the updating of $P_M$ depends on such choice. Nevertheless some 
maintenance aspects of $P_M$ can be mentioned. The update of $P_M$ should be 
delayed as much as possible due to the fact that updates on $P_M$ force an 
update on $B$.
A great advantage on the modularity of $P$ is that it allows among other
things, to change preconditioning strategy ($P_M$) mid-way between two updates.

Preliminary experimentation shows that most of the big changes for 
$\big\|M_{\ell}-M_{\ell-1}\big\|_1$ occur at the beginning and specially between
two external iterations due to the Lagrangian multipliers and the external 
penalty parameters being updated. Given the convergence of ALM and supposing the 
problem has solution, iterates will converge asymptotically
to a point where updating $P_M$ will no longer be necessary.
The reported experiments suggest to use a lax threshold for the update of $P_M$.
This choice promotes frequent updates only 
at the beginning of the resolution while avoiding unnecessary updates towards 
the end.

\subsubsection{On the update of matrix $B$}\label{sec:updB}
Matrix $B$ is associated with the constraint Jacobian and as such its form is
described by the active non-relaxed constraints and the rank-1 correction 
artifacts.
Given the recursive nature of the assembly of $B$, it is possible to establish 
predictive update strategies that reduce costs. The idea is the following.
Let $V_{m-2}$ be the second to last column of $V$ and the gradient associated 
with the inequality constraint  
$\c_{m-2}(z_{\ell})$ and also suppose that it has a very small infeasibility measure. Most
probably the next iterate will inactivate this constraint forcing the discarding
of $V_{m-2}$. This induces and update \emph{only} to the last two columns of $B$.
Now, if $V_{m-2}$ where to be the first column of $V$ then the induced 
update would affect \emph{all} the columns of $B$ having a much greater cost.
This observation suggest orderings over the columns of $V$ that 
potentially reduce the cost of updating $B$. A practical ordering could be induced by the
infeasibility measure of each constraint 
\[ \big|c_i(z)\big|, \;\; i\in\set{E}, \qquad\text{y}\qquad 
\max\!\big(0,\; c_i(z)\big), \;\; i\in\set{I}. \]
The spirit behind this ordering is to leave for the end of the recursion those 
columns whose  associated constraints will (possibly) soon be discarded.
% This ordering as the effect of lowering the cost for discarding those elements
% from $B$.
This ordering could be further improved by taking into account the norm of each
associated gradient.
It is also very convenient to leave at the end of the recursion the $q$ columns
associated with BFGS-type corrections. These elements can have a vivid transit 
state since initially iterations may or may not fulfill the condition 
$s^{\kern -1pt \T}y<10^{-8}\|s\|\,\|y\|.$
Although, it has been observed that close to a solution the BFGS correction can
always be applied.
Numerical experimentation shows that as iterates approach a solution, changes in
$V$ diminish down to a point where updating $B$ is not necessary and can be 
recycled successfully.

\subsection{Comments}
Accelerating Cauchy-type methods require to consider preconditioning matrices as
approximations to the inverse of the Hessian. 
% Under this context, 
Acceleration
enriches the descent direction $-\nabla\F(x)$ with second order information
from the approximation $H$ of $\nabla^2\F(x)$. In others words, the enriched
direction is obtained by solving 
\[H_k\,d_k=-\nabla\F(x_k),\]
but this is exactly the same task as applying acceleration under Newton-like
choices.
Hence the instructions on how to apply and when to update the (linear) 
preconditioner are analogous for the nonlinear case.
Unfortunately the enriched (preconditioned) direction is not always a 
descent direction. In these cases, the enriched direction is discarded in 
favor of $-\nabla\F(x_k)$, albeit, the work invested in building the 
approximation $P\approx H^{-1}$ should not be discarded: it could potentially be 
recycled the next time an enriched direction is to be computed.

The assembly and maintenance of matrix $B$ loses appeal and stops being 
attractive in the presence of a large number of active constraints, say $m\approx n$ or even $m>n$. 
In this scenario, we have $H=M+W$ where $W$ is dense and possibly rank
complete or near complete.
As iterates start closing-in to a solution, the number of active constraints
should decrease down to a point where the use of $B$ is practical. 
This suggest to handle acceleration for problems with a large amount of 
constraints  in a two-stage approach. 
% The following briefly mentions a heuristic.
% \begin{enumerate}
% \item 

A possible heuristic is to use Quasi-Newton directions induced by $V_{\set{K}}(z)$ where $\set{K}$
is the set of constraints indexes for the $0<K\ll{}n$ elements with greatest
infeasibility measure. The quadratic model for the direction is
\[ V_K{V_K}^{\!\!\!\!\T}\; d = -\nabla\F(z),\quad \text{ or even }\quad %\]
% or possibly enriched with the Lagrangian Hessian
% \[ 
\big[M + V_K{V_K}^{\!\!\!\!\T}\big]\; d = -\nabla\F(z). \]
Under these two choices, it is possible to use the proposed scheme. 
It is important to note that before finding the final search space, the set
$\set{K}$ will be changing inducing unfavorable frequent updates on matrix $B$.
% 
%  \item 
% A second choice is to set $P=P_M$, using only $M$ ($\nabla^2\!\!\lagrangian(z)$).
% The direction is found by solving
% $$P_M^{-1}H\;d = -P_M^{-1}\nabla\F(z).$$ 
% As soon as the number of active non-relaxable constraints drops to a 
% reasonable amount, then $P$ is updated to account for $\rho\,VV^{\T}$.
% \end{enumerate}

It is important to consider alternative techniques that do not use the 
explicit form of the constraint Jacobian but can tackle with the dense 
Gauss-Newton matrix $W$. This topic is considered open for future study.

\input{experimentacion_short}

\section{Concluding remarks}
An acceleration scheme for the Augmented Lagrangian method was presented. 
The associated preconditioner ($P$) exploits the explicit form of the Augmented
Lagrangian Hessian ($H$) without estimating its inverse. The strategy is modular
and uses two main ingredients. An auxiliary preconditioner ($P_M$) associated 
with the Lagrangian Hessian ($M$) and a storage matrix ($B$) related to the constraints 
Jacobian matrix and possible rank-one corrections of $H$. The preconditioner takes
inspiration in the Sherman-Morrison identity and Miller's inverse formula.
A virtue of this scheme is that the acceleration strategy is agnostic to the class
of preconditioner used for $P_M$. For special choices of approximations to the 
Lagrangian Hessian, the preconditioner  is the exact inverse.
Quasi-Newton and Secant-type approximations are encouraged, associated 
low-rank BFGS corrections are absorbed in the $B$ matrix and do not require 
explicit handling. The quality of $P$ is determined by the quality of $P_M$ and 
the relaxation over the constraints during the assembly of $B$. Some 
characteristics of $P$ are induced by $P_M$.
The explicit handling of the constraints is attractive when $m<n$. This is not
necessarily the case when $m\approx n$ or $m>n$, since $VV^{\T}$ is possibly 
rank complete and the paradigm of using the explicit form of all constraints is 
no longer practical. Other methods that take advantage of the aggregated form of 
the Gauss-Newton matrix are recommended.
From the implementation point of view, the preconditioner along with the 
updating strategies can be directly used on Newton-type as well as on 
Projected Gradient methods. This reduces the coding effort and allows code 
recycling.

Numerical experiments gave insight on the quality of the preconditioner and 
general behavior of the scheme inside ALM. Initial results reveal that when $P_M$
is poor, so is $P$. Optimal assembly parameters are problem specific and not 
universal.
Constraint relaxation is quite delicate and required careful handling. 
Arbitrary relaxation proved to be a poor choice. Recommended relaxation 
strategies take into account infeasibility measure and gradient size.
Further study on the ordering of constraints during the assembly of $B$ is
recommended.
During the experimentation a wide range of incomplete direct and inverse 
factorizations where used for the auxiliary preconditioner $P_M$. Results 
show expected behavior and highlight the agnostic quality of the scheme.
Good preconditioner update strategies monitor the changes between to iterates. 
Refreshing on a fixed number of iterations is not recommended.
The acceleration scheme was successfully coupled to two inner solvers:
a Truncated-Newton machinery for unconstrained sub-problems and
a projected gradient-type for convex constrained problems.
From Tables~\ref{tbl:e5}, \ref{tbl:e6} and \ref{tbl:e7} 
it can be observed that inside the ALM context, preconditioning is not only 
attractive, but can also be the only alternative that produces a solution.
Also can be inferred that attractive assembly parameter values, relaxation and 
update tolerances for Quasi-Newton machinery are also good for the Spectral
Projected Gradient case.

The natural next step is the implementation of the scheme in a low level 
language and incorporation inside an existing ALM implementation with the 
objective of further experimenting on larger problems.

\begin{multicols}{2}
\section*{References}\scriptsize
\bibliographystyle{elsarticle-num}
\bibliography{biblio}
\end{multicols}

\end{document}

%% file: alg-MLA.tex
\begin{algorithm}
\caption{Augmented Lagrange Method~\cite{libroMLA}.}
\label{alg:MLA}
\begin{algorithmic}
%  \REQUIRE
\STATE \textbf{Input: }
 $\lambda_{\min}<\lambda_{\max}, \mu_{\max}>0, \gamma>1, 0<\tau<1,
 \bar{\lambda}_1\in[\lambda_{\min}, \lambda_{\max}]^m, \bar{\mu}_1\in[0,\mu_{\max}]^p$,
 $\rho_1>0$. $k\leftarrow1$.
 \STATE\textbf{Step~1} Find approximate solution $x_k\in\set{R}^n$ of~\eqref{eq:minlagaug}.
 \STATE\textbf{Step~2} Stop if $x_k$ satisfies \eqref{eq:KKT1}--\eqref{eq:KKT3}. % is solution.
 \STATE\textbf{Step~3} Update Lagrange multipliers: 
 \[\lambda_{k+1} = \bar{\lambda}_k+\rho_k\h(x_k) \quad \text{and} \quad 
   \mu_{k+1} = \Big(\bar{\mu}_k+\rho_k\g(x_k)\Big)_+\]
 \STATE\textbf{Step~4} Update $\rho$:
 \STATE Let $V_{i,k} = \min\big(-g_i(x_k), \; \bar{\mu}_{i,k}/\rho_k\big),\quad i=1:p$. 
 \IF{$k=1 \;\vee\; \max\Big\{\|\h(x_k)\|,\|V_{k}\|\Big\}\le\tau\;\max\Big\{\|\h(x_{k-1})\|,\|V_{k-1}\|\Big\}$}
 \STATE $\rho_{k+1}\ge\rho_k$.
 \ELSE
 \STATE $\rho_{k+1}=\gamma\rho_k$.
 \ENDIF
 \STATE\textbf{Step~5} Safeguard multipliers: $\bar{\lambda}_{k+1}\in[\lambda_{\min}, \lambda_{\max}]^m, \bar{\mu}_{k+1}\in[0,\mu_{\max}]^p$.
 \STATE\textbf{Step~6} Increment iteration counter $k\leftarrow k+1$ and go to Step~1.
\end{algorithmic}
\end{algorithm}

%% file: experimentacion_short.tex
\section{Numerical Experiments}
All experiments were run using Matlab\textregistered\ R2012a on an 
Intel\textregistered\ Core\texttrademark\ i7-2640M CPU 
@ 2.80GHz with 8~GB of memory.
We start by examining the quality of $P$ as a preconditioner for $H$ and
its efficiency at solving linear systems of the form $Hx=y$. 
Some experiments involving constraint relaxation and preconditioner update
strategies follow. We finalize by solving unconstrained and box-constrained
problems from the CUTEst\cite{Gould2015cutest} data-set.

\subsection{Spectral Properties of the Preconditioned Matrix}\label{sec:eigen}
\begin{figure}[htbp]
{\centering\begin{tabular}{p{0.22\linewidth}p{0.22\linewidth}p{0.22\linewidth}p{0.22\linewidth}}
\multicolumn{2}{c}{\footnotesize Matrix of size 100 with 10 constraints} &
\multicolumn{2}{c}{\footnotesize Matrix of size 100 with 50 constraints}\\
\multicolumn{2}{c}{\includegraphics[width=0.44\linewidth]{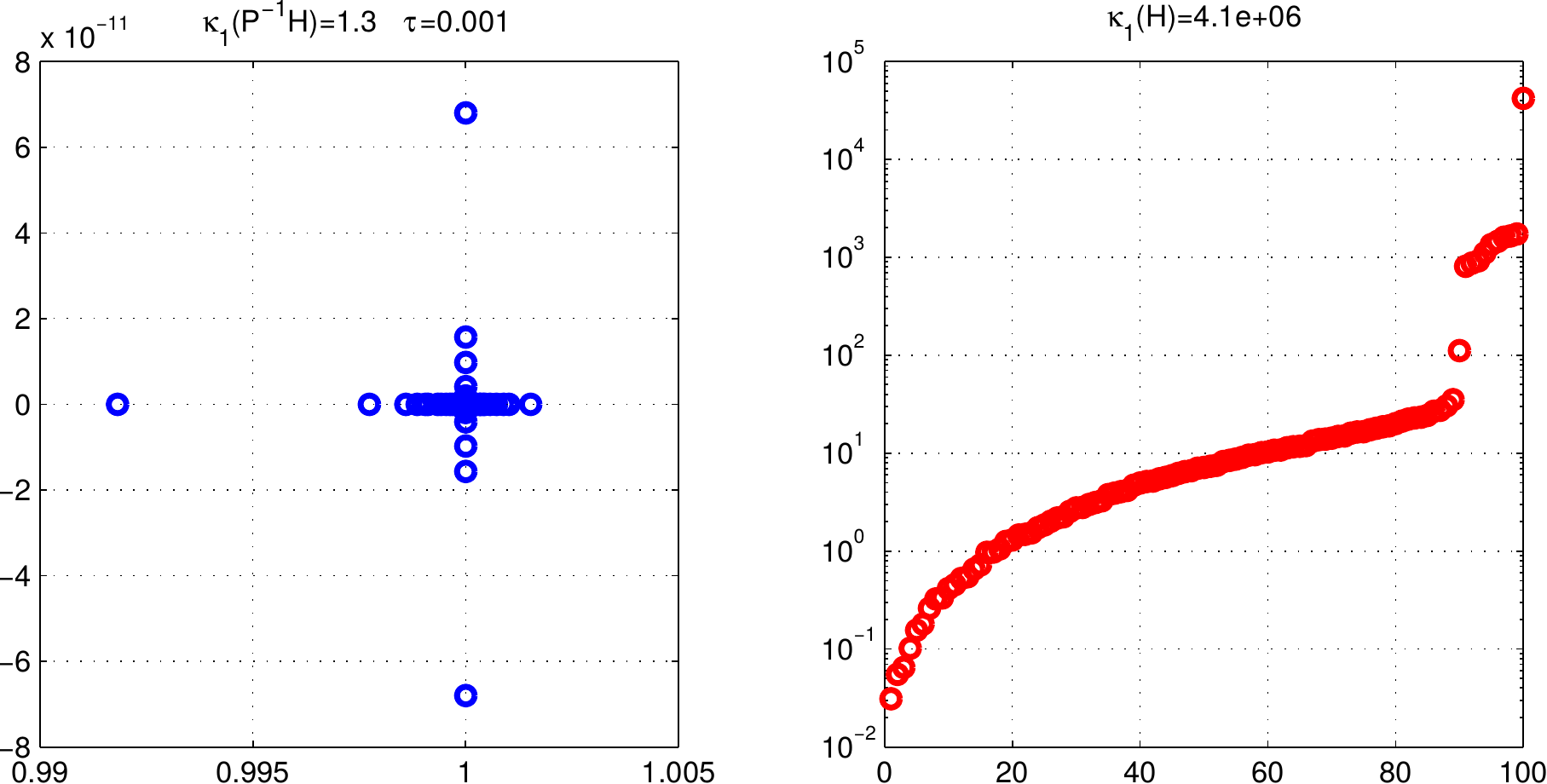}} &
\multicolumn{2}{c}{\includegraphics[width=0.44\linewidth]{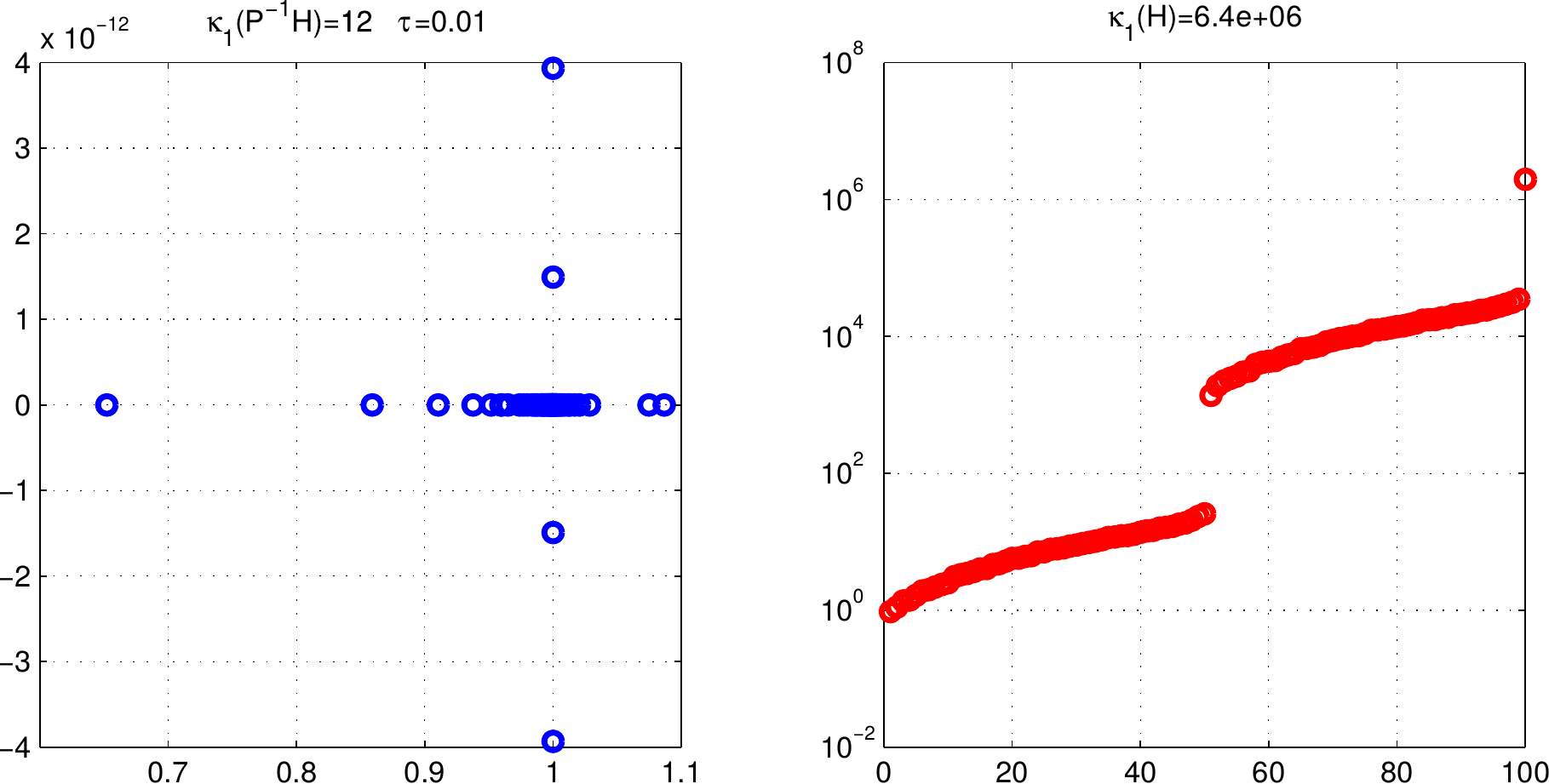}}\\
\centering \footnotesize (a) & 
\centering \footnotesize (b) & 
\centering \footnotesize (c) & 
\centering \footnotesize (d)
\end{tabular}}
\caption{Spectral distribution of $P^{-1}H$ and $H$.
Figures (a) and (c), 
show in blue color the eigenvalues of $P^{-1}H$ 
grouped near the identity $(1, 0)$.
Figures (b) y (d) show in red color the eigenvalue distribution of $H$ in logarithmic scale.}
\label{fig:eig1}
\end{figure}

We wish to understand the spectral properties of the preconditioned matrix
 $P^{-1}H=P^{-1}\big[M+\rho\,VV^{\T}\big]$. Quality metric is based on condition 
 number and spectra of  $P^{-1}H$. For these experiments, the auxiliary 
 preconditioner $P_M$ is of the family of Robust Incomplete Factorization of 
 type SAINV~\cite{benzi1996sparse,benzi2000robust,benzi2003robust,
benzi2003approximate} and is considered a black-box that executes the 
matrix-vector product $M^{-1}Y$. 
Table~\ref{tbl:e1} and Figure~\ref{fig:eig1}
report two particular experiments on sparse random matrices of size 100 with
$m\in\{10,50\}$ random constraints. These experiments use a diverse range of 
values for the external penalty $\rho$ and dropping $\tau$ parameters.
\input{tbl-e1}From the table and figure it can be observed that the higher the quality of
$P_M$ (smaller $\tau$) the spectrum of $P^{-1}H$ accumulates around the 
identity, albeit, if $P_M$ is poor then the conditions of $P^{-1}H$ and $H$
are comparable. Curiously for these problems, higher values of $\rho$ seem to
have a favorable effect over the condition of $P^{-1}H$. Increased values of
$\rho$ worsen the condition of $H$ while promoting the one of $P^{-1}H$.

The next set of experiments involve solving the linear system of equations
$Hx=y$ with $H=M+\rho\,VV^{\T}$ using the Conjugate Gradients method. Tested 
matrices are sparse random and from the Matrix~Market~\cite{MatrixMarket} 
collection. All matrices are forced to be real SPD. Constraints are random with 
N($0, 1$) distribution. Convergence tolerance is set to $10^{-8}$.
Obtained result are reported in Table~\ref{tbl:e2}.
Dropping parameters $\tau_1$ y $\tau_2$ are associated with $P_M$. 
Column labeled with ``$m$'' represents the number of constraints and the columns
associated with $\nnz(Z)$ give an idea of the density of $P_M$. Columns ``CG'' 
and ``PCG'' report the number of iterations required.
From the table an evident influence of 
$\rho$ over the conditioning of $H$ is observed. 
The higher the quality of $P_M$, the more evident is the acceleration.
On some problems, for very large values $\rho$ the number of CG iterations is
surprisingly low and preconditioning stops being practical. Conversely, in some
cases where $H$ is very ill conditioned, preconditioning not only is very 
effective but is the only alternative that converges. In PCG context, 
increment in $\rho$ generally implies reduction in iteration count.

\input{tbl-e2}

\subsection{Constraint Relaxation \& Update}
Updating the preconditioner on each iteration is prohibitively expensive and even
unnecessary. The convergence of ALM guarantees that sub-problems will tend to 
be similar up to the point where recycling the preconditioner is possible. 
For inequality constraints, as soon as the iterate is feasible or strongly feasible,
the number of active constraints drops drastically. This leads to propose 
cheap update strategies that discard elements from the $B$ matrix.
Active constraint relaxation must be done with certain care. Preliminary 
experiments suggest it is not easy to establish the \emph{contribution} of each
active constraint to the Gauss-Newton matrix.
With the idea of understanding the influence of each constraint (column of $V$)
on the quality of the preconditioner $P$, spectral properties of $P^{-1}H$ are
monitored while assembling $P$ with a sub-set of columns of $V$. For these
experiments $P_M$ can be considered $M^{-1}$.
Preconditioner $P$ is assembled by increasing the amount of columns of $V$ 
included. The order of inclusion is given by the norm of each column,
$\|v_1\|\ge\|v_2\|\ge\cdots\ge\|v_m\|>0$.
Figure~\ref{fig:ccs} shows the condition number of $P^{-1}H$ in function of
$\rho$ and the number of columns of $V$ used in the assembly of $P$.
Poor behavior can be observed, in order to maintain good quality most of the
components of $V$ must be used. Attractive relaxation strategies have to account 
at least for the norm and infeasibility measure of the components of $V$. 
This suggests to assemble $P$ using only those active constraints where
$\left\|\sqrt{\rho_k}\,\nabla\c_i(x_{\ell})\right\|\ge\varepsilon_v$ and 
$\left|\c_i(x_{\ell})\right|\ge\varepsilon_c$.

\input{fig-ccs}

Updating strategies for $P$ must take advantage of the convergence of ALM.
At first changes to $M$ and $V$ are assumed to be \emph{big} but will tend
to smooth as the sub-problems start to be similar.
In order to understand the associated costs of diverse updating strategies for 
$P$, the following experiments are conducted. Let $\delta_M$ and $\delta_v$ be
updating tolerances for changes in matrices $M$ and $V$ between two
consecutive iterations. We solve problem C4~\cite{thesis-phd-usb-sajo-2016} 
varying parameters 
$\delta_M$, $\delta_v$, $\varepsilon_v$ (relaxation on $\nabla\c(x)$) and $\tau$
(dropping tolerance for $P_M$), while observing the amount of updates required 
to produce a solution.
Table~\ref{tbl:e4} reports a sequence of experiments on problem C4 to understand
updating strategies/cost relation. The cost of each configuration is given by
 CG/MinRes iterations and the amount of updates on $M$ and $V$
required to find a solution. These two counters are antagonistic. In general terms,
low CG/MinRes iteration count corresponds to a high update count for $M$ and $V$. The 
idea is to find a compromise between these two. Reported numerical results
in Table~\ref{tbl:e4} give an intuitive overview of attractive configurations. 
% % Results favor the use of lax update tolerances.
% 
\input{tbl-e4}

Table~\ref{tbl:C4-RBUWX} shows in detail the dynamics of certain parameters
while solving problem C4 using QN machinery. The problem is rigged to have
ten variables with nine inequality and one equality constraints.
Descent direction is found using MinRes. The auxiliary preconditioner 
$P_M$ is ILU of type Crout with dropping tolerance $10^{-7}$.
Column labeled ``Update type'' specifies the three reasons for updating $P$: 
``M'' and ``V'' indicate that changes between two consecutive $M$ or $V$ 
matrices where greater that $\delta_M$ and/or $\delta_v$ respectively, while 
``V*'' indicates  a \emph{forced} update to $H$ due to BFGS-style corrections 
($|\set{V}|$ changes from 10 to 12). It can also indicate that the correction 
cannot be applyed in the current iteration ($|\set{V}|$ changes from 
12 to 10). $|\set{V}|$ indicates the quantity of columns of matrix $V$ in use 
and rnnz($Y$) indicates the density of $P_M$ relative to $M$.
From Table~\ref{tbl:C4-RBUWX} it can be observed that initially the changes 
between two consecutive matrices $M$ and $V$ are substantial, hence the update
type ``MV''. Further in the resolution of the problem, changes to $M$ start to 
attenuate then diminish for $V$. In general terms, this behavior ca be
considered characteristic. In particular, for problem C4 after iteration 4 of
Lagrangian iteration 3 it was no longer necessary to update $P$ which 
represents a 44.8\% reduction in unnecessary updates.
% 
\input{tbl-C4-RBUWX}
Obtained results show the convenience of updating $P$ via $\delta_M$ and 
$\delta_v$ tolerances over strategies that update on predefined iterations
e.~g. every Lagrangian iteration.
Unfortunately, experimentation showed that attractive update tolerances is problem dependent 
and requires individual adjustment. For some problems these can be lax and still
give good acceleration results, specially for changes between elements of $V$.
It is also attractive to relax the \emph{small} columns of $V$ since  
it potentially reduces costs and can efficiently be controlled using a single
additional relaxation parameter. Adding and removing columns from $V$ is a frequent
task. This requires cheap machineries that handle efficiently the 
addition and removal of elements within the matrix $B$. Furthermore, the assembly
order inside of $B$ can take advantage of a prediction ingredient based on 
 constraint gradient size and infeasibility measure.

\subsection{Unconstrained Sub-problems}
In each Lagrangian iteration, the problem to solve is
\[ \text{minimize} \F(x) = \lagrangiana_{\rho_k}(x, \lambda_k)\; \text{subject to } x\in\set{R}^n,\; \text{with $\rho_k$ and $\lambda_k$ fixed.}\]
Two type of quadratic models where used for the direction. Newton (NW) $H_{\ell} = \nabla^2\F(z_{\ell})$ and
Quasi-Newton (QN)
$$H_{\ell}=\nabla^2\f(x_{\ell})  + \sum_{i\in\set{A}}\nabla\c_i(x_{\ell})\nabla\c_i(x_{\ell})^{\T} + \sigma I + \text{BFGS},$$
where $\set{A}$ is the set of active non-relaxed constraints indexes, $\sigma$ the inverse Rayleigh quotient and BFGS represents rank-one corrections.
The auxiliary preconditioner $P_M$ used is from the ILU family. 
Table~\ref{tbl:e5}
reports the obtained experimental results. Acceleration is evident when
contrasting columns labeled ``Itpd'' 
and ``Itd''. For this set of problems, the preconditioning scheme achieved 
acceleration factors between 1.3 and 9.6.
\input{tbl-e5}
On some problems the QN direction did not converge.
Most of the computational effort is invested in updating $P_M$, this promotes
approximations to $M$ that are easy invertible, 
such as band-diagonal approximations.
Preconditioners used in the NW and QN models are also used in the 
Spectral Gradient method (SG)\footnote{SPG code adapted from the
TANGO project~\cite{birgin2000nonmonotone,birgin2001algorithm}.} and its
preconditioned variant (PSG).
Comparing the column ``Itin'' of the rows ``SG'' and ``PSG'' for each problem 
in Table~\ref{tbl:e6} acceleration in iteration count is evident, but if
comparing computational cost, preconditioning is not always attractive.

\input{tbl-e6}

\subsection{Box-constrained Sub-problems}
For these experiments, the problem to solve is
\[ \min \F(x) = \lagrangiana_{\rho_k}(x, \lambda_k)\quad \text{subject to }    
\ell_i\le x \le u_i,\quad \text{with $\rho_k$ and $\lambda_k$ fixed.}\]
In order to understand the efficiency of the scheme, the sub-problems are
solved using SPG and its preconditioned variant PSPG.
Table~\ref{tbl:e7} reports numerical experiments,
similar results to Table~\ref{tbl:e6} were obtained. Acceleration is clear when
applying the preconditioner (HS105, HS111 y HS112).
In some problems using the preconditioner reduces iteration count but 
increments overall cost (EXTRASIM, HS41 y HS63). This last aspect indicates that
fine-tuning is required  in order to avoid premature activation of the
preconditioner.
Premature activation of $P$ does not increase total iteration count but does 
have a negative impact of overall computational cost.
The idea is to find a compromise between updating $P$ too frequently
and iterating a large amount of times.
A conservative strategy would be to delay the activation of the 
preconditioner up until the iterates are \emph{close} to a solution. Another
alternative could be to activate preconditioning as soon as the iterates start
to be \emph{very feasible} thus promoting the acceptance of the 
preconditioned direction.

\input{tbl-e7}

%% file: tbl-e1.tex
\begin{table}[hbtp]
\caption{Experimental results for solving the linear system $Hx=y$, 
with $H=M+\rho\,VV^{\T}$ for random $M\in\set{R}^{100\times100}$ and 
$V\in\set{R}^{100\times m}$.}\label{tbl:e1}
\begin{multicols}{2}\centering

\begin{tabular}{p{2.5ex}p{9ex}p{6ex}p{11ex}p{11ex}}
\toprule
$m$ & $\rho$ & $\tau$ & $\kappa_1(H)$ & $\kappa_1(P^{-1}H)$  \\\midrule
 10 & 1.5 & 0.1    & $2.4\times10^{5}$ & $5.4\times10^{4}$ \\
 10 & 15.5 & 0.1   & $5.8\times10^{5}$ & $4.9\times10^{4}$ \\
 10 & 154.8 & 0.1  & $4.1\times10^{6}$ & $4.9\times10^{4}$ \\
 10 & 1548.3 & 0.1 & $4.0\times10^{7}$ & $4.9\times10^{4}$ \\
 10 & 15483  & 0.1& $4.0\times10^{8}$ & $4.9\times10^{4}$ \\
\midrule
 10 & 1.5 & 0.001    & $2.4\times10^{5}$ & 1.3 \\
 10 & 15.5 & 0.001   & $5.8\times10^{5}$ & 1.3 \\
 10 & 154.8 & 0.001  & $4.1\times10^{6}$ & 1.3 \\
 10 & 1548.3 & 0.001 & $4.0\times10^{7}$ & 1.3 \\
 10 & 15483  & 0.001& $4.0\times10^{8}$ & 1.3 \\
 \bottomrule
\end{tabular}

\begin{tabular}{p{2.5ex}p{9ex}p{6ex}p{11ex}p{11ex}}
\toprule
$m$ & $\rho$ & $\tau$ & $\kappa_1(H)$ & $\kappa_1(P^{-1}H)$  \\\midrule
 50 & 1.5 & 0.1    & $2.0\times10^{4}$ & $5.5\times10^{3}$ \\
 50 & 15.5 & 0.1   & $7.3\times10^{4}$ & $3.4\times10^{3}$ \\
 50 & 154.8 & 0.1  & $6.5\times10^{5}$ & $3.6\times10^{3}$ \\
 50 & 1548.3 & 0.1 & $6.4\times10^{6}$ & $3.7\times10^{3}$ \\
 50 & 15483  & 0.1& $6.4\times10^{7}$ & $3.7\times10^{3}$ \\
 \midrule
 50 & 1.5 & 0.01    & $2.0\times10^{4}$ & 20 \\
 50 & 15.5 & 0.01   & $7.3\times10^{4}$ & 13 \\
 50 & 154.8 & 0.01  & $6.5\times10^{5}$ & 12 \\
 50 & 1548.3 & 0.01 & $6.4\times10^{6}$ & 12 \\
 50 & 15483  & 0.01& $6.4\times10^{7}$ & 12 \\
\bottomrule
\end{tabular}
\end{multicols}
\end{table}

%% file: tbl-e2.tex
\begin{table}\noindent\small
\caption{Experimental results for solving the linear system $Hx=y$ using
Conjugate Gradients method.}\label{tbl:e2}
\begin{tabular}{lccccccclllcc}
\toprule
Name &$n$&$m$&$\tau_1$&$\tau_2$&$\nnz(Z)$&$\frac{\nnz(Z)}{n^2}$&$\frac{\nnz(Z)}{\nnz(M)}$&$\rho$&$\kappa(H)$&$\kappa(P^{-1}H)$&CG&PCG\\\midrule
Sparse&1000&1&0.05&0.075&2040&0.002&0.674&1&12681509&74245&n/c&22\\
% Sparse&1000&1&0.05&0.075&2040&0.002&0.674&10&126038198&74251&n/c&20\\
Sparse&1000&1&0.05&0.075&2040&0.002&0.674&100&1259615534&74251&172&14\\
Sparse&1000&1&0.05&0.075&2040&0.002&0.674&1000&12595389913&74251&58&9\\
Sparse&1000&100&0.1&0.1&1656&0.002&0.684&100&153329943&12&56&2\\
Sparse&1000&50&0.1&0.1&1666&0.002&0.684&1&2650446&122&173&9\\
Sparse&1000&50&0.1&0.1&1666&0.002&0.684&100&264040355&121&54&3\\
Sparse&2000&1&0.05&0.075&2865&0.001&0.746&1&70457448&70219&n/c&15\\
Sparse&2000&1&0.05&0.075&2865&0.001&0.746&1000&70285577870&70214&37&3\\
bcspwr01&49&20&0.1&0.1&131&0.086&1&1&208&1.3&21&5\\
% bcspwr01&49&20&0.1&0.1&131&0.086&1&100&18446&1.2&34&4\\
bcspwr01&49&20&0.1&0.1 &131&0.086&1&10000&1841063&1.2&27&2\\
bcspwr01&49&1&0.18&0.18&142&0.059&0.85&1&171&15&31&17\\
% bcspwr01&49&1&0.18&0.18&142&0.059&0.85&10&1379&15&30&16\\
bcspwr01&49&1&0.18&0.18&142&0.059&0.85&100&13474&15&30&14\\
% bcspwr01&49&1&0.18&0.18&142&0.059&0.85&1000&134420&15&23&13\\
bcspwr01&49&1&0.18&0.18&142&0.059&0.85&10000&1343884&15&21&9\\
bcspwr01&49&1&0.18&0.18&142&0.059&0.85&100000&13438526&15&9&6\\
bcspwr02&49&1&0.18&0.18&159&0.066&0.952&1&600&53&33&18\\
% bcspwr02&49&1&0.18&0.18&159&0.066&0.952&10&4688&52&32&17\\
bcspwr02&49&1&0.18&0.18&159&0.066&0.952&100&45618&51&32&15\\
% bcspwr02&49&1&0.18&0.18&159&0.066&0.952&1000&454922&51&31&14\\
bcspwr02&49&1&0.18&0.18&159&0.066&0.952&10000&4547969&51&20&12\\
% bcspwr02&49&1&0.18&0.18&159&0.066&0.952&100000&45478435&51&10&9\\
bcspwr02&49&1&0.18&0.18&159&0.066&0.952&1000000&454783100&51&5&3\\
bcspwr02&49&10&0.1&0.1&157&0.065&0.94&1&294&1.6&25&6\\
% bcspwr02&49&10&0.1&0.1&157&0.065&0.94&100&22598&1.4&33&4\\
bcspwr02&49&10&0.1&0.1&157&0.065&0.94&10000&2253720&1.4&27&2\\
bcspwr02&49&30&0.1&0.1&157&0.065&0.94&1&355&1.3&21&5\\
bcspwr02&49&30&0.1&0.1&157&0.065&0.94&10&2720&1.2&39&4\\
bcsstm03&112&1&0.7&0.7&112&0.009&1&1&3003.3&1&42&1\\
bcsstm03&112&1&0.7&0.7&112&0.009&1&1000000&23186836&1&19&1\\
bcsstm03&112&1&0.7&0.7&112&0.009&1&10000000&231867364&1&7&1\\\bottomrule
\end{tabular}
\end{table}

%% file: fig-ccs.tex
\begin{figure}[htbp]
\centering
\begin{tabular}{cc}
\includegraphics[width=0.473\linewidth]{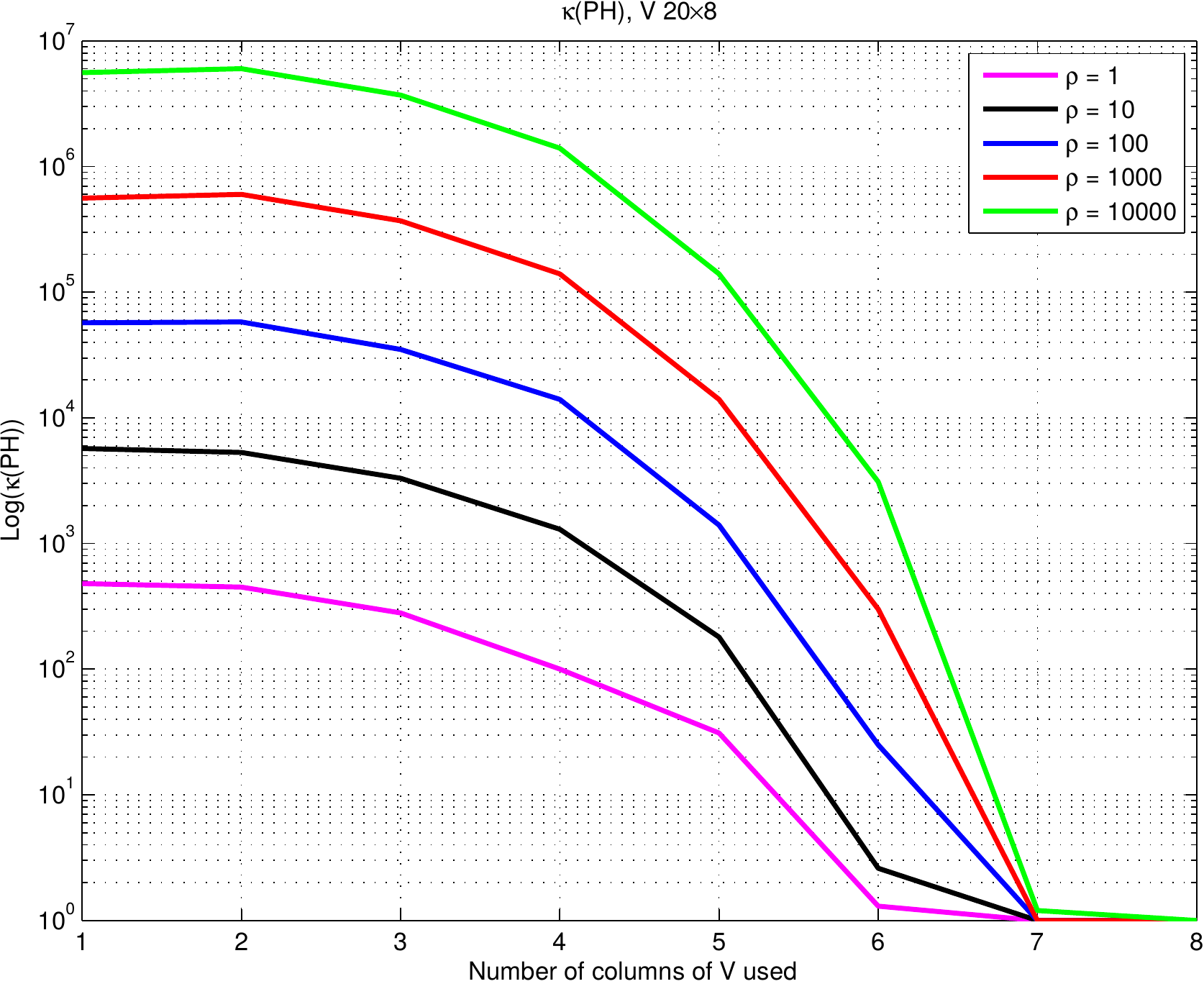} &
\includegraphics[width=0.473\linewidth]{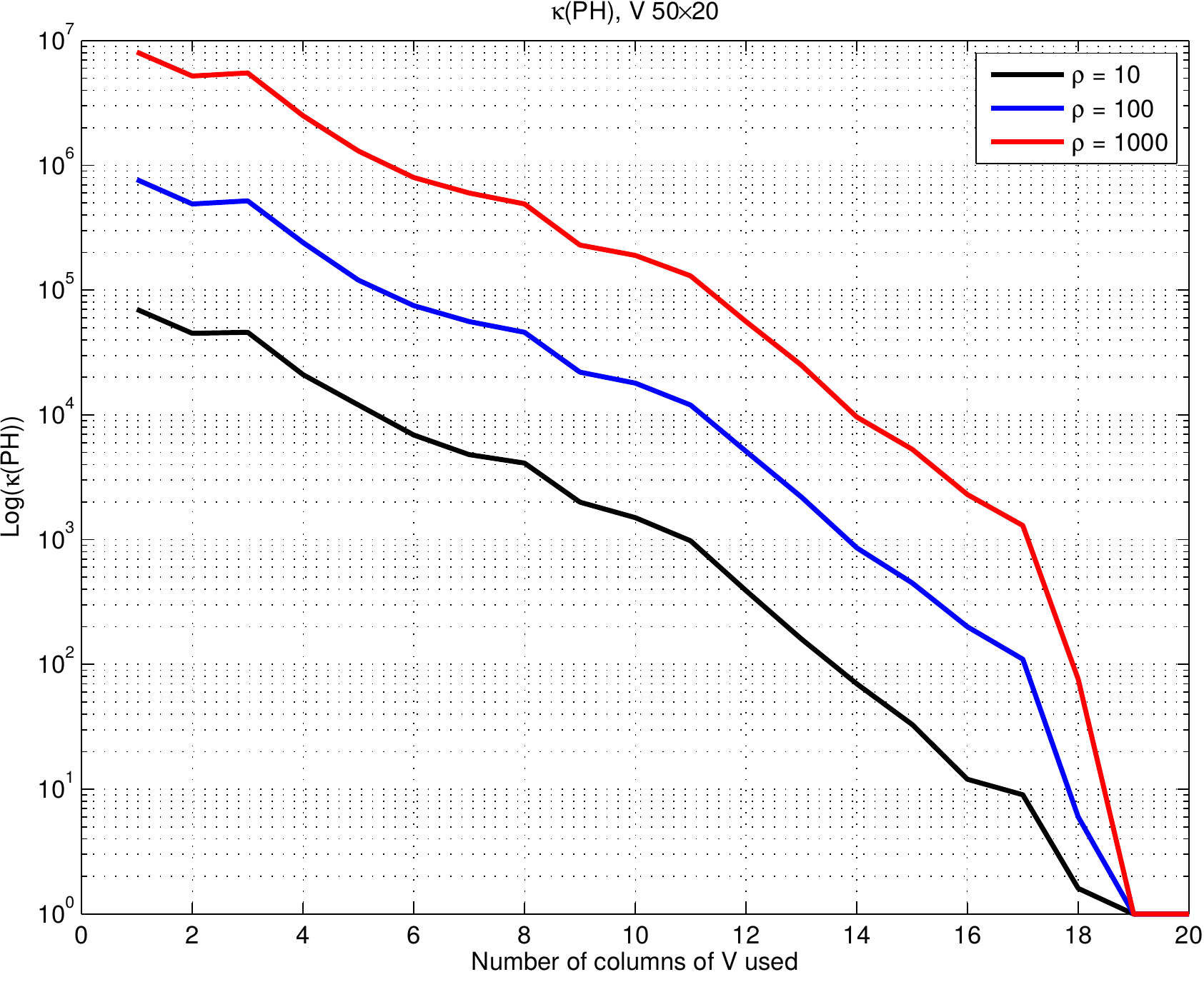}\\
\footnotesize (a) & \footnotesize  (b)
\end{tabular}
\caption{Accumulative influence of the columns of $V$ in the quality of the
preconditioner $P$.
In (a) the problem is of size 20 with 8 active constraints. For $P$ to be 
competitive, the first 9 columns of $V$ are required in the assembly of $B$.
In (b) the problem has dimension 50 with 20 constraints. In this case the first 
19 columns of $V$ are needed.
\label{fig:ccs}}
\end{figure}

%% file: tbl-e4.tex
\begin{table}[tbhp]
 \caption{Numerical experiments to understant iteration costs based on update 
 strategies and modulation of the parameters  $\delta_M$, $\delta_v$, 
 $\varepsilon_v$ and $\tau$ for problem C4. 
 \label{tbl:e4}}
\centering {\small
\begin{tabular}{lccccccccccc}\toprule
           &&&&&&&&\multicolumn{3}{c}{Updates}\\
Strategy & $\delta_M$      & $\delta_v$      & $\varepsilon_v$ & $\tau$          & Ext. it & Int. it & GC/MR it  & Total & $M$  & $V$ \\
\midrule
Auto        & $1\cdot10^{-1}$ & $1\cdot10^{-1}$ & $1\cdot10^{-2}$ & $10^{-6}$       & 37    & 98    & 317       & 56        & 20     & 36\\  % SUMMARY
Auto        & $1\cdot10^{-1}$ & $1\cdot10^{-2}$ & $1\cdot10^{-2}$ & $10^{-6}$       & 37    & 98    & 317       & 62        & 20     & 42\\  % SUMMARY
Auto        & $1\cdot10^{-1}$ & $1\cdot10^{-2}$ & $1\cdot10^{-3}$ & $10^{-6}$       & 37    & 98    & 317       & 62        & 20     & 42\\  % SUMMARY
Auto        & $1\cdot10^{-1}$ & $1\cdot10^{-2}$ & $1\cdot10^{-3}$ & $10^{-9}$       & 37    & 98    & 317       & 62        & 20     & 42\\  % SUMMARY
Auto        & $1\cdot10^{-2}$ & $1\cdot10^{-2}$ & $1\cdot10^{-3}$ & $10^{-6}$       & 37    & 98    & 258       & 65        & 30     & 35\\  % SUMMARY
Auto        & $2\cdot10^{-1}$ & $1\cdot10^{-2}$ & $1\cdot10^{-3}$ & $10^{-9}$       & 37    & 98    & 345       & 62        & 18     & 44\\  % SUMMARY
Auto        & $2\cdot10^{-1}$ & $1\cdot10^{-2}$ & $5\cdot10^{-3}$ & $10^{-9}$       & 37    & 98    & 345       & 62        & 18     & 44\\  % SUMMARY
Auto        & $5\cdot10^{-1}$ & $1\cdot10^{-2}$ & $5\cdot10^{-3}$ & $10^{-6}$       & 37    & 98    & 368       & 62        & 17     & 45\\  % SUMMARY
Auto        & $5\cdot10^{-1}$ & $5\cdot10^{-1}$ & $1\cdot10^{-1}$ & $10^{-3}$       & 37    & 98    & 399       & 27        & 17     & 10\\  % SUMMARY
Auto        & $5\cdot10^{-1}$ & $5\cdot10^{-1}$ & $1\cdot10^{-2}$ & $10^{-6}$       & 37    & 98    & 394       & 27        & 17     & 10\\  % SUMMARY
Auto        & $5\cdot10^{-1}$ & $5\cdot10^{-3}$ & $1\cdot10^{-3}$ & $10^{-3}$       & 37    & 98    & 373       & 65        & 17     & 48\\  % SUMMARY
Auto        & $5\cdot10^{-2}$ & $1\cdot10^{-2}$ & $1\cdot10^{-3}$ & $10^{-9}$       & 37    & 98    & 284       & 62        & 22     & 40\\  % SUMMARY ** \rowcolor{highlightcolor}
\multicolumn{5}{l}{Update every external iteration}                                 & 37    & 98    & 533       & 35        & 3      & 32\\  % SUMMARY
\multicolumn{5}{l}{Assemble $P$ only once}                                          & 37    & 98    & 1156      & 1         & 1      & 0\\   % SUMMARY
\multicolumn{5}{l}{Newton / direct method}                                          & 37    & 98    & ---       & ---       & ---    & ---\\ % SUMMARY
\bottomrule
\end{tabular}}
\end{table}

%% file: tbl-C4-RBUWX.tex
% Version: 8.1
\begin{table}[hbtp]
\caption{\label{tbl:C4-RBUWX}%
Evolution of the iterates while solving problem C4 using ALM using QN-type directions.}
 \centering \begin{tabular}{crcccccrr}
\toprule
 L. it & N. it  & \multicolumn{1}{p{4em}}{\centering PMR it (MR it)} & \multicolumn{1}{p{3em}}{\centering Update type} & $|\mathbb{V}|$ & $\rho$ & rnnz($Y$) & $\|M_k-M_{k-1}\|_1$ & $\|v_k-v_{k-1}\|_1$\\
\midrule
 1 &  1 &    1 (10)  &  MV  & 10 & 100 & 1.6 &\multicolumn{1}{c}{---}&\multicolumn{1}{c}{---}  \\
 1 &  2 &    1 (7)   &  MV*  & 12 & 100 & 1.6 & $1\times10^{   3}$    & $1.16\times10^{3}$ \\
 1 &  3 &    1 (5)   &  MV  & 12 & 100 & 1.6 & $13.9$                & $788$ \\
 1 &  4 &    1 (5)   &  MV  & 12 & 100 & 1.6 & $2.95$                & $850$ \\
 1 &  5 &    1 (5)   &  MV  & 12 & 100 & 1.6 & $9.42\times10^{-1}$   & $318$ \\
 1 &  6 &    1 (5)   &  MV  & 12 & 100 & 1.6 & $5.37\times10^{-1}$   & $354$ \\
 1 &  7 &    1 (5)   &  MV  & 12 & 100 & 1.6 & $7.13\times10^{-1}$   & $129$ \\
 1 &  8 &    1 (5)   &  MV  & 12 & 100 & 1.6 & $1.56\times10^{-1}$   & $189$ \\
 1 &  9 &    1 (7)   &  MV  & 12 & 100 & 1.6 & $978$                 & $50.3$ \\
 1 & 10 &    1 (8)   &  MV  & 12 & 100 & 1.6 & $9.81$                & $2.1$ \\
 \midrule
 2 &  1 &    1 (10)  &  MV*  & 10 & 100 & 1.6 & $12.6$                & $0.0$ \\
 2 &  2 &    1 (10)  &  MV*  & 12 & 100 & 1.6 & $2.59$                & $4.91$ \\
 \midrule
 3 &  1 &    1 (10)  &  MV*  & 10 & 100 & 1.6 & $2.53$                & $0.0$ \\
 3 &  2 &    1 (10)  &  M   & 10 & 100 & 1.6 & $1.38\times10^{-1}$   & $0.0$ \\
 3 &  3 &    3 (10)  &  V*  & 12 & 100 & 1.6 & $3.56\times10^{-3}$   & $3.76$ \\
 3 &  4 &    3 (10)  &  V   & 12 & 100 & 1.6 & $2.22\times10^{-3}$   & $3.77\times10^{-1}$ \\
 3 &  5 &    3 (10)  & ---  & 12 & 100 & 1.6 & $1.11\times10^{-3}$   & $9.78\times10^{-3}$ \\
 3 &  6 &    3 (10)  & ---  & 12 & 100 & 1.6 & $5.53\times10^{-4}$   & $6.36\times10^{-3}$ \\
 3 &  7 &    3 (10)  & ---  & 12 & 100 & 1.6 & $2.76\times10^{-4}$   & $3.20\times10^{-3}$ \\
 3 &  8 &    3 (10)  & ---  & 12 & 100 & 1.6 & $1.38\times10^{-4}$   & $1.60\times10^{-3}$ \\
 3 &  9 &    3 (10)  & ---  & 12 & 100 & 1.6 & $6.90\times10^{-5}$   & $8.02\times10^{-4}$ \\
 3 & 10 &    3 (10)  & ---  & 12 & 100 & 1.6 & $3.45\times10^{-5}$   & $4.01\times10^{-4}$ \\
 3 & 11 &    3 (10)  & ---  & 12 & 100 & 1.6 & $1.72\times10^{-5}$   & $2.01\times10^{-4}$ \\
 3 & 12 &    3 (10)  & ---  & 12 & 100 & 1.6 & $8.62\times10^{-6}$   & $1.00\times10^{-4}$ \\
 3 & 13 &    3 (10)  & ---  & 12 & 100 & 1.6 & $4.31\times10^{-6}$   & $5.02\times10^{-5}$ \\
 3 & 14 &    2 (10)  & ---  & 12 & 100 & 1.6 & $2.16\times10^{-6}$   & $2.51\times10^{-5}$ \\
 3 & 15 &    2 (10)  & ---  & 12 & 100 & 1.6 & $1.08\times10^{-6}$   & $1.25\times10^{-5}$ \\
 3 & 16 &    2 (10)  & ---  & 12 & 100 & 1.6 & $5.39\times10^{-7}$   & $6.27\times10^{-6}$ \\
 3 & 17 &    2 (10)  & ---  & 12 & 100 & 1.6 & $5.39\times10^{-7}$   & $2.89\times10^{-6}$ \\\bottomrule
\end{tabular}
\end{table}

%% file: tbl-e5.tex
\begin{table}[htbp]
\centering
\caption{\label{tbl:e5}%
Numerical results for solving unconstrained problems using
ALM coupled with Newton-type methods. Columns labeled by 
``ItL.'', ``Itin'', ``Itpd'' and ``Itd'' 
indicate the number of external, internal, MinRes and accelerated MinRes 
iteration count. The sum of columns ``AcM'' and ``AcV'' reports the times
the preconditioner $P$ was updated while each column shows individual 
contribution.}

% \small
\begin{tabular}{lcccccccccc}\toprule
Name & $n$ & $m$  & Method & ItL. & Itin & Itpd & Itd & AcM & AcV & Time (s) \\\midrule %& Identificación\\\midrule
C4-1 & 10 & 10       & NW  & 11 & 52 & 82 & 420 & 35 & 1 & 2.4 \\% & HWKEV \\       -SMMILU
C4-1 & 10 & 10       & NW  & 10 & 53 & 83 & 432 & 34 & 2 & 2.5 \\% & WWBRG \\       -SMMILU
C4-1 & 10 & 10       & QN  & 45 & 143 & 311 & 1365 & 40 & 19 & 5.7 \\% & KNWKX \\   -SMMILU
C4-1 & 10 & 10       & QN  & 45 & 143 & 176 & 1365 & 107 & 6 & 5.9 \\% & HQPMQ \\   -SMMILU
C4-R-10& 10 & 10     & QN  & 3 & 29 & 55 & 252 & 14 & 2 & 1.0 \\% & RBUWX\\         -SMMILU
C4-R-20& 20 & 10     & QN  & 3 & 48 & 265 & 617 & 33 & 6 & 1.6 \\% & PBIKF\\        -SMMILU
C4-R-25& 25 & 10     & QN  & 4 & 55 & 112 & 900 & 37 & 18 & 1.7 \\% & XBXEQ\\       -SMMILU
\midrule
BT3 & 5 & 3          & NW & 8 & 84 & 84 & 807 & 1 & 3 & 1.3 \\% & SGMQV \\          -SMMILU
BT3 & 5 & 3          & NW & 7 & 86 & 86 & 539 & 1 & 1 & 1.4 \\% & BBCTW \\          -SMMILU
BT3 & 5 & 3          & QN & 8 & 129 & 182 & 651 & 21 & 58 & 1.8 \\% & LDSCH \\      -SMMILU
BT3 & 5 & 3          & QN & 8 & 129 & 193 & 651 & 20 & 50 & 1.7 \\% & LSJGA \\      -SMMILU
\midrule             
BT8 & 5 & 2          & NW & 5 & 72 & 87 & 148 & 27 & 0 & 1.1  \\% & QNHDO \\        -SMMILU
BT8 & 5 & 2          & QN & 5 & 88 & 108 & 169 & 26 & 8 & 1.3 \\% & XEGJB \\        -SMMILU
% \midrule             
% BT9 & 4 & 2          & NW & 11 & 159 & 293 & 381 & 37 & 4 & 2.3 \\% & GSEGC \\      -SMMILU
% BT9 & 4 & 2          & QN & n/c &---&---&---&---&---&---\\
\midrule
BT11 & 5 & 3         & NW & 10 & 140 & 322 & 707 & 38 & 4 & 2.1 \\% & SFKQI \\      -SMMILU
BT11 & 5 & 3         & NW & 8 & 114 & 134 & 680 & 112 & 0 & 1.8 \\% & LELDB \\      -SMMILU
% BT11 & 5 & 3         & QN & n/c&---&---&---&---&---&---\\
% \midrule             
% BT12 & 5 & 3         & NW & 6 & 54 & 125 & 411 & 12 & 0 & 0.9  \\% & ROFPC \\       -SMMILU
% BT12 & 5 & 3         & QN & 7 & 146 & 281 & 587 & 7 & 48 & 1.8 \\% & RABQO \\       -SMMILU
% \midrule             
% HS48 & 5 & 2         & NW & 2 & 31 & 31 & 155 & 1 & 0 & 0.6 \\% & RFCOK \\          -SMMILU
% HS48 & 5 & 2         & QN & n/c &---&---&---&---&---&---\\ 
% \midrule             
% MAKELA4 & 21 & 40    & NW & 37 & 228 & 228 & 786 & 1 & 27 & 6.1 \\% & SRJPE \\      -SMMILU
% MAKELA4 & 21 & 40    & QN & n/c &---&---&---&---&---&---\\
\bottomrule
% \multicolumn{11 }{l}{{\footnotesize Problemas marcados con un asterisco (*) indican que están re-escalados.}}
\end{tabular}

\end{table}

%% file: tbl-e6.tex
\begin{table}[htbp]
\centering
\caption{\label{tbl:e6}%
Numerical results for solving unconstrained problems using
ALM coupled with SG and it preconditioned variant PSG. 
Column labels are described in the caption of 
Table~\ref{tbl:e5}.}
% \small
\begin{tabular}{lcccccccccc}\toprule
Name & $n$ & $m$ & Method 
                        & ItL. &Itin& AcM & AcV & Time (s) \\\midrule% & Identificación\\\midrule
C4-1 & 10 & 10 &  SG & n/c& ---   & ---   & ---  & ---\\% & UDFUT\\         
C4-1 & 10 & 10 & PSG & 1 & 21 & 3 & 0 & 0.3  \\% & OGORF \\          -SMMILU
% C4-1 & 10 & 10 & PSG & 1 & 38 & 20 & 1 & 0.6 \\% & KKTCD \\          -SMMILU
% C4-1 & 10 & 10 & PSG & 1 & 38 & 18 & 2 & 0.6 \\% & CIELI \\          -SMMILU
\midrule
BT3 & 5 & 3 &   SG  & 9 & 486 & --- & --- & 1.6 \\% & KJSTE \\              
% BT3 & 5 & 3 &  PSG & 9 & 36  & 8 & 0 & 0.5  \\% & LKPRS \\          SMMILU 
% BT3 & 5 & 3 &  PSG & 9 & 68  & 8 & 27 & 0.8 \\% & GQPDC \\          SMMILU 
BT3 & 5 & 3 &  PSG & 7 & 19  & 12 & 0 & 0.4 \\% & QKKOB \\          SMMILU 
\midrule
BT8 & 5 & 2 &   SG & 3 & 29  & --- & --- & 0.3 \\% & DJDSU \\              
% BT8 & 5 & 2 &  PSG & 3 & 25  & 2 & 1 & 0.4 \\% & WWKFS \\           SMMILU 
BT8 & 5 & 2 &  PSG & 3 & 17  & 5 & 0 & 0.3 \\% & DPVMQ \\           SMMILU 
\midrule
BT9 & 4 & 2 &   SG & 9 & 193 & --- & --- & 0.8 \\% & HQPMQ \\              
BT9 & 4 & 2 &  PSG & 9 & 32  & 25 & 0 & 0.5  \\% & PEDXE \\         SMMILU 
% BT9 & 4 & 2 &  PSG & 9 & 111 & 68 & 14 & 1.4 \\% & AKHDE \\         SMMILU 
\midrule
BT11 & 5 & 3 &  SG & 7 & 426 & --- & --- & 1.5 \\% & TVXAU \\              
BT11 & 5 & 3 & PSG & 7 & 28 & 28 & 0 & 0.4    \\% & FLVNU \\       -SMMILU
\midrule
BT12 & 5 & 3 &  SG & 7 &1222& --- & --- & 4.3 \\% & CLUUG \\              
% BT12 & 5 & 3 & PSG & 6 & 14 & 14 & 0 & 0.3 \\% & FNPKE \\          -SMMILU
BT12 & 5 & 3 & PSG & 6 & 16 & 12 & 0 & 0.3 \\% & WBCDD \\          -SMMILU
% BT12 & 5 & 3 & PSG & 6 & 33 & 14 & 4 & 0.5 \\% & ONBWR \\          -SMMILU
% BT12 & 5 & 3 & PSG & 6 & 31 & 10 & 12 & 0.5 \\% & RBUWX \\         -SMMILU
\midrule
HS48 & 5 & 2 &  SG  & 2 & 131& --- & --- & 0.5 \\% & RJFGS \\              
% HS48 & 5 & 2 & PSG & 1 & 2  & 2 & 0 & 0.2 \\% & BQORX \\           -SMMILU
HS48 & 5 & 2 & PSG & 1 & 2  & 2 & 0 & 0.2 \\% & JGVTL \\           -SMMILU
\midrule
MAKELA4 & 21 & 40 &  SG & 2 & 5 & 0 & 0 & 0.1 \\% & TQHWA \\
MAKELA4 & 21 & 40 & PSG & 2 & 3 & 3 & 0 & 0.3 \\% & QAGBC \\
\bottomrule
\end{tabular}
\end{table}

%% file: tbl-e7.tex
\begin{table}[htbp]
\centering
\caption{Numerical results for solving convex constrained problems using
ALM coupled with SPG and PSPG. Column labels are described in the caption of 
Table~\ref{tbl:e5}.\label{tbl:e7}}

% \small
\begin{tabular}{lcccccccccc}\toprule
Name & $n$ & $m$ & Method 
                         &ItL.& Itin& AcM &AcV& Time (s) \\\midrule% & Identificación\\\midrule
AIRPORT & 84 & 42 &  SPG & n/c& --- & --- &---&---   \\% & \\\midrule
AIRPORT & 84 & 42 & PSPG & 81 & 304 & 159 & 0 & 18.9 \\% & TQHWA \\              -SMMILU 
AIRPORT & 84 & 42 & PSPG & 81 & 304 & 159 & 0 & 22.9 \\% & KJSTE \\              -SMMILU 
\midrule
EXTRASIM & 2 & 1  &  SPG & 1  & 16  & --- &---& 0.2 \\% & XNDDG \\\midrule  
EXTRASIM & 2 & 1  & PSPG & 1  & 2   & 2   & 0 & 0.2 \\% & UGTFW \\              -SMMILU             
\midrule
HS41 & 4 & 1      &  SPG & 4  & 78  & --- &---& 0.4 \\% & DNLAI \\                  
HS41 & 4 & 1      & PSPG & 4  & 78  & 6   & 0 & 0.7 \\% & KVEGD \\\midrule  -SMMILU 
\midrule
HS63 & 3 & 2      &  SPG & 4  & 125 &---  &---& 0.5 \\% & LTHTN \\                  
HS63 & 3 & 2      & PSPG & 4  & 26  & 16  & 0 & 0.5 \\% & WWRJA \\\midrule  -SMMILU 
\midrule
HS90 & 4 & 1      & SPG  & 1  & 10  &---  &---& 0.1 \\% & NHREQ \\                  
HS90 & 4 & 1      & PSPG & 1  & 8   & 1   & 0 & 0.2 \\% & USMEJ \\\midrule  -SMMILU 
\midrule
HS105 & 8 & 1     &  SPG & 2  & 505 & --- &---& 4.5 \\% & DPUXN \\                   
HS105 & 8 & 1     & PSPG & 1  & 100 & 82  & 0 & 2.1 \\% & MGWVJ \\           -SMMILU 
HS105 & 8 & 1     & PSPG & 1  & 80  & 49  & 0 & 1.7 \\% & AAFKD \\\midrule   -SMMILU 
\midrule
HS111 & 10 & 3    &  SPG & 12 & 5999& --- &---& 21.7\\% & JPAWL \\                    
HS111 & 10 & 3    & PSPG & 12 & 56  & 41  & 0 & 0.8 \\% & EVGUT \\            -SMMILU 
HS111 & 10 & 3    & PSPG & 12 & 86  & 36  & 0 & 1.2 \\% & CAGWQ \\            -SMMILU 
\midrule
HS112 & 10 & 3    &  SPG & 12 & 1567& --- &---& 5.7 \\% & PHHEC \\                    
HS112 & 10 & 3    & PSPG & 12 & 37  & 35  & 0 & 0.8 \\% & RTOCN \\\midrule    -SMMILU 
\midrule
LOOTSMA & 3 & 2   &  SPG & 1  & 34  & --- &---& 0.6 \\% & LHOVQ \\                     
LOOTSMA & 3 & 2   & PSPG & 1  & 34  & 1   & 0 & 0.5 \\% & RNDJJ \\\midrule     -SMMILU 
\bottomrule
\end{tabular}
\end{table}